\documentclass[11pt,epsfig]{article}
\usepackage{amssymb}
\usepackage{amsmath}
\usepackage{amsthm}
\usepackage{epsfig}
\usepackage{hyperref}
\usepackage{amsfonts}
\usepackage{color}
\usepackage[margin=1.2in]{geometry}
\usepackage{lipsum,titletoc}

\titlecontents{section} 
[2.5em]                 
{\rmfamily}
{\contentslabel{2.3em}} 
{\hspace*{-2.3em}} {\titlerule*[1pc]{.}\contentspage}

\titlecontents{subsection} 
[4.8em]                 
{\rmfamily}             
{\contentslabel{2.3em}} 
{\hspace*{-2.3em}} {\titlerule*[1pc]{.}\contentspage}

\newtheoremstyle{thmm}{1.5ex plus 1ex minus .2ex}{1.5ex plus 1ex minus
.2ex}{\rmfamily}{}{\bfseries}{}{1em}{} \theoremstyle{thmm}
\newtheorem{theorem}{Theorem}[section]
\newtheorem{lemma}{Lemma}[section]
\newtheorem{proposition}{Proposition}[section]

\renewenvironment{proof}[1][Proof]{\noindent\textit{#1. } }{\hfill$\square$}

\allowdisplaybreaks

\newcommand{\nn}{\nonumber}

\def\d{\delta}

\def\R{\mathbb{R}}

\def\d{{\rm d}}

\begin{document}

\title{\bf
Uniform BMO estimate of parabolic equations and global well-posedness of the  thermistor problem
\\[10pt]}

\setcounter{footnote}{0}

\author{
Buyang  Li \footnote{Department of Mathematics,
Nanjing University, Nanjing, P.R. China.
\newline \indent~ {\tt Email
\!address:~buyangli@nju.edu.cn} } ~ and ~Chaoxia Yang }

\date{}

\maketitle

\begin{abstract}
Global well-posedness of the time-dependent (degenerate) thermistor
problem remains open for many years. In this paper, we solve the
problem by establishing a uniform-in-time BMO estimate of
inhomogeneous parabolic equations. Applying this estimate to the
temperature equation, we derive a BMO bound of the temperature
uniform with respect to time, which implies that the electric
conductivity is a $A_2$ weight. The H\"{o}lder continuity of the
electric potential is then proved by applying the De
Giorgi--Nash--Moser estimate for degenerate elliptic equations with
$A_2$ coefficient. Uniqueness of solution is proved based on the
established regularity of the weak solution. Our results also imply
the existence of a global classical solution when the initial
and boundary data are smooth. \\[5pt]

\end{abstract}\medskip

\tableofcontents

\section{Introduction}
\setcounter{equation}{0} The thermistor problem refers to the
heating of a conductor, with temperature-sensitive electric
conductivity, by electric current.
Let $\phi$ be the electric potential and let ${\bf E}=\nabla\phi$ be
the electric field. The electric current ${\bf J}$ is related to the
electric field via ${\bf J}=\sigma(u){\bf E}$, where
$\sigma(u)$ is the electric conductivity of the conductor, dependent upon the temperature $u$.
The heat produced (per unit volume) by the electric current is given
by Joule's law: ${\bf E}\cdot{\bf J}=\sigma(u)|\nabla\phi|^2$, and the
conservation of charge is described by $\nabla\cdot{\bf J}=0$.

Let $\Omega$ denote the domain possessed by the conductor. Based on
the above formulations, the temperature $u$ and the electric
potential $\phi$ are governed by the equations
\begin{align}
&\frac{\partial u}{\partial t}-\nabla\cdot(\kappa(u)\nabla
u)=\sigma(u)|\nabla\phi|^2, \label{e-heat-1}
\\[3pt]
&-\nabla\cdot(\sigma(u)\nabla\phi)=0, \label{e-heat-2}
\end{align}
for $x\in\Omega$ and $t>0$, where $\kappa(u)$ is the thermal
conductivity. In this paper, we
consider the above equations with the Dirichlet boundary/initial
conditions:
\begin{align}
\label{BC}
\begin{array}{ll}
u(x,t)=g(x,t),\quad \phi(x,t)=h(x,t)~~
&\mbox{for}~~x\in\partial\Omega~~\mbox{and}~~t>0,\\[3pt]
u(x,0)=u_0(x)~~ &\mbox{for}~~x\in\Omega .
\end{array}
\end{align}

The mathematical expressions of $\sigma(u)$ and $\kappa(u)$ depend
on the materials. For some semiconductors, the electric resistivity
$\rho(u)=1/\sigma(u)$ can be approximately expressed as
\cite{Macklen}
$$
\rho(u)=\sigma_0e^{q/u}u,
$$
and the thermal conductivity $\kappa(u)$ can be regarded as constant
(independent of $u$). For metallic conductors, the electric
conductivity and the thermal conductivity obey the Wiedemann--Franz
law \cite{Sommerfeld}:
\begin{align*}
\frac{\kappa(u)}{\sigma(u)}=Lu ,
\end{align*}
where $L=2.44\times10^{-8}{\rm W\Omega K^{-2}}$ is the Lorentz
number. In general, the electric resistivity of metals increases as
temperature grows. At high temperatures, the electric resistivity
increases approximately linearly with temperature:
$$
\rho(u)=\rho_0[1+\alpha(u-u_R)],
$$
where $u_R$ is some reference temperature and $\alpha$ is called the
temperature coefficient of resistivity. If the temperature does not vary
much, the above linear formula is often used. More precisely, the
electric resistivity is give by the Bloch--Gr\"{u}neisen formula
\cite{Ziman}:
$$
\rho(u)=\rho(0)+A\biggl(\frac{u}{\Theta}\biggl)^n\int_0^{\frac{\Theta}{u}}\frac{s^n}{(e^s-1)(1-e^{-s})}\d
s ,
$$
where $A$, $\Theta$ and $n\geq2$ are all positive physical
constants.

For both metals and semiconductors, the electric conductivity
$\sigma(u)$ tends to zero as the temperature $u$ grows to infinity.
The elliptic equation (\ref{e-heat-2}) is thus possibly degenerate,
which leads to severe difficulties for the analysis of the
coupled system.

The non-degenerate assumption
$\sigma_1\leq\sigma(u)\leq \sigma_2$ is often used to
simplify the problem.  Mathematical analysis for such
non-degenerate problem has been
studied by many authors in the last two decades.
Existence of weak solutions was studied by Antontsev and Chipot
\cite{AC}, Allegretto and Xie \cite{AX} and Cimatti \cite{Cimatti2}. 
With the same non-degenerate assumption, Elliott and Larsson
\cite{EL} proved the existence of strong solutions for the 2D
problem by using the energy method (and uniqueness follows). The 3D problem
is much more difficult. To deal with the 3D problem, one has to fully
explore and make use of the coupling of the equations. The milestone
was acheived by Yuan and Liu \cite{YL,Yuan2}, who proved the existence of
$C^\alpha$ solutions for the 3D problem by using the method of Layer
potentials. Yin \cite{Yin} obtained the same result by using the techniques
of Campanato spaces.
Their results imply the existence of classical
solutions when the boundary and initial data are smooth.

Without the non-degenerate assumption, the problem becomes much more
difficult. Xu \cite{Xu} proved partial regularity of the solution,
i.e. the solution is smooth in an open subset $D\subset\Omega$ whose
complement $\Omega\backslash D$ is a set of measure zero. Later Xu
\cite{Xu1} proved existence of solutions with bounded temperature
when the boundary potential is small enough, i.e.
$\|h\|_{L^\infty(\partial\Omega\times(0,T))}$ is small enough.
Hachimi and Ammi \cite{HA} proved existence of weak solutions by the
monotonicity-compacity method. Montesinos and Gallego \cite{MG1,MG2}
proved existence of ``capacity solutions'' by considering a new
formulation with the transformation $\Phi=\sigma(u)\nabla\phi$.
Uniqueness of the weak solution and existence of global classical
solutions remain open. Overall, the main difficulty of the
degenerate problem is the lack of a $L^\infty$ bound for the
temperature $u$.

In this paper, we overcome this difficulty by establishing a
uniform-in-time ${\rm BMO}$ estimate for inhomogeneous parabolic
equations with possibly discontinuous coefficients. Applying this
estimate to the temperature equation, we obtain a uniform-in-time
${\rm BMO}$ bound of the temperature $u$, as a substitute of the
$L^\infty$ bound. Based on the BMO bound of the temperature, we
further prove that the electric conductivity $\sigma(u)$ is a $A_2$
weight uniform in time. The H\"{o}lder continuity of the electric
potential $\phi$ is then proved by applying the De
Giorgi--Nash--Moser estimate for degenerate elliptic equations with
$A_2$ coefficient. The H\"{o}lder continuity of the temperature is
proved by using the H\"{o}lder continuity of the electric potential.
Existence of a weak solution in a bounded Lipschitz domain is
proved, and uniqueness of the weak solution is proved based on the
established regularity of the solution. Our results also imply the
existence of a global classical solution when the initial and
boundary data are smooth.

For interested readers, we refer to \cite{AL, AY, EL, LS, LS2, ZW} for
numerical methods and numerical analysis of the thermistor problem.

The rest part of this paper is organized in the following way. In
Section \ref{notations} we introduce the notations to be used in
this paper and in Section \ref{dfjksl7879} we present our main
results. In Section \ref{fjkesluio99}, we establish a
uniform-in-time BMO estimate for the solutions of inhomogeneous
parabolic equations, and in Section \ref{sec5} we present H\"{o}lder
estimates of parabolic equations in terms of the Campanato spaces.
Based on the estimates obtained in Section \ref{fjkesluio99} and
Section \ref{sec5}, we prove global existence and uniqueness of a
weak solution to the degenerate thermistor problem in Section
\ref{fsdjkljl78687}. Conclusions are drawn in Section 7.

\section{Notations}\label{notations}
\setcounter{equation}{0}

Before we present our main results, we define the notations to be used in this paper.

Let $n$ be a fixed positive
integer and let $B_R(x_0)$ denote the ball of radius $R$ centered at
the point $x_0\in\R^n$. Let $\Omega$ be a bounded Lipschitz domain
in $\R^n$, i.e. $\Omega$ is a bounded domain in $\R^n$ and for any
$y\in\partial\Omega$, there exists a ball $B_R(y)$ such that through
a rotation of coordinates (if necessary),
$$
B_R(y)\cap\Omega=\{(x_1,\cdots,x_n)\in
B_R(y):x_n>\varphi(x_1,\cdots,x_{n-1})\} ,
$$
where $\varphi:\R^{n-1}\rightarrow \R$ is a Lipschitz continuous
function. For a bounded Lipschitz domain, there exists a positive
constant $R_\Omega$ and a finite number of balls
$B_{R_\Omega}(y_1)$, $B_{R_\Omega}(y_2)$, $\cdots$,
$B_{R_\Omega}(y_m)$ such that $\partial\Omega\subset
\cup_{j=1}^mB_{R_\Omega/2}(y_j)$ and through a rotation of
coordinates (if necessary),
$$
B_{2R_\Omega}(y_j)\cap\Omega=\{(x_1,\cdots,x_n)\in
B_{2R_\Omega}(y_j):x_n>\varphi_j(x_1,\cdots,x_{n-1})\}
$$
for some Lipschitz continuous function
$\varphi_j:\R^{n-1}\rightarrow \R$.

For any integer $m\geq 0$, $1\leq p\leq\infty$ and $0<\alpha<1$, let
$W^{m,p}(\Omega)$ and $C^{m+\alpha}(\overline\Omega)$ denote the
usual Sobolev space and H\"{o}lder space \cite{Adams}, respectively,
and let $C^{m+\alpha}(\Omega)$ denote the space of functions which
belong to $C^{m+\alpha}(\overline B)$ for any closed ball $\overline
B\subset\Omega$. Let $C^{m+\alpha}_0(\overline\Omega)$ be the
subspace of $C^{m+\alpha}(\overline\Omega)$ consisting of functions
vanishing on the boundary $\partial\Omega$.

Let $|D|$ denote the Lebesgue measure for any measurable subset $D$ of $\R^n$, and let $B_R(x_0)$ denote the ball of radius $R$ centered at
the point $x_0\in\R^n$. Let $\Omega$ be a bounded Lipschitz domain
in $\R^n$. We say that a positive locally integrable function $w$ defined on $\R^n$ is a $A_2$ weight if
$$
\sup_{B\subset\R^n}\biggl(\frac{1}{|B|}\int_Bw(x)\d x\biggl)\biggl(\frac{1}{|B|}\int_B\frac{1}{w(x)}\d x\biggl) \leq C
$$
for some positive constant $C$, where the supremum extends over all balls in $B$ in $\R^n$.

For any measurable
subset $D$ of $\R^n$, we let
$f_{D}=\frac{1}{|D|}\int_{D}f(x)\d x$ denote the average of $f$ over
$D$. For $1\leq p<\infty$ and $0\leq \theta\leq 1$, let
$L^{p,\theta}(\Omega)$ denote the Morrey space of measurable functions $f$ such that

\begin{align*}
\|f\|_{L^{p,\theta}(\Omega)}:=\sup_{B_R(x_0)}\biggl(\frac{1}{R^{n\theta}}\int_{B_R(x_0)\cap\Omega }|f(x)|^p\d
x\biggl)^\frac{1}{p}<\infty,
\end{align*}
where the supremum above extends over all balls $B_R(x_0)$ with
$x_0\in\overline\Omega$ and $0<R<R_\Omega$. For $1\leq p<\infty$ and
$1\leq\theta<\infty$, let ${\cal L}^{p,\theta}(\Omega)$ denote the
Campanato space of functions bounded (or vanishing for $\theta>1$)
on the boundary $\partial\Omega$, equipped with the norm
\begin{align*}
\|f\|_{{\cal L}^{p,\theta}(\Omega)}:
&= \sup_{
B_R(x_0)\cap\Omega}\biggl(\frac{1}{R^{n\theta}}\int_{
B_R(x_0)\cap\Omega}|f(x)|^p\d
x\biggl)^\frac{1}{p}\\
&~~~+\sup_{B_R(y_0)\cap\Omega}\biggl(\frac{1}{R^{n\theta}}\int_{B_R(y_0)\cap\Omega}|f(x)-f_{B_R(y_0)\cap\Omega}|^p\d
x\biggl)^\frac{1}{p},
\end{align*}
where the supremum above extends over all balls with
$x_0\in\partial\Omega$, $y_0\in\Omega$ and $0<R<R_\Omega$, and we set
$\overline{\rm BMO}={\cal L}^{1,1}(\Omega)$.

For any fixed $T>0$, we set $\Omega_T=\Omega\times(0,T]$ and $\Gamma_T=\partial\Omega\times(0,T]$. For any
point $(x_0,t_0)\in\R^{n+1}$, we set
$Q_R(x_0,t_0)=B_R(x_0)\times(t_0-R^2,t_0]$ as the parabolic cylinder centered
at $(x_0,t_0)$ of radius $R$. For integers $m,n\geq 0$,
$0<\alpha,\beta<1$ and any open subset $Q\subset \Omega_T$, let
$C^{m+\alpha,n+\beta}(\overline Q)$ denote the anistropic H\"{o}lder
space of functions, equipped with the norm
$\|f\|_{C^{m+\alpha,n+\beta}(\overline
Q)}:=\|f\|_{L^\infty(Q)}+|f|_{C^{m+\alpha,n+\beta}(\overline Q)}$,
where
\begin{align*}
|f|_{C^{m+\alpha,n+\beta}(\overline Q)}&=\sum_{|\gamma|=
m}\sup_{\begin{subarray}{l}
(x,t)\in Q\\
(y,s)\in Q
\end{subarray}}\frac{|D^\gamma_x f(x,t)-D^\gamma_x
f(y,s)|}{|x-y|^\alpha +|t-s|^\beta}\\
&~+\sum_{|\gamma|= n}\sup_{\begin{subarray}{l}
(x,t)\in Q\\
(y,s)\in Q
\end{subarray}}\frac{|D^\gamma_t f(x,t)-D^\gamma_t
f(y,s)|}{|x-y|^\alpha+|t-s|^\beta},
\end{align*}
and set
$C^\alpha(\overline\Omega_T)=C^{\alpha,\alpha}(\overline\Omega_T)$.
Let $C^{m+\alpha,n+\beta}_0(\overline Q)$ denote the subspace of
$C^{m+\alpha,n+\beta}(\overline Q)$ with functions vanishing on the
boundary $\partial\Omega$. Let $C^\infty(\overline Q)$ denote the
space of functions whose partial derivatives up to all orders are
uniformly continuous on $\overline Q$. Let
$C^{m+\alpha,n+\beta}(\Omega_T)$ and $C^\infty(\Omega_T)$ denote the
space of functions which are in $C^{m+\alpha,n+\beta}(\overline Q)$
and $C^\infty(\overline Q)$ for any closed cylinder
$\overline Q\subset\Omega_T$, respectively. For any measurable
subset $Q$ of $\R^{n+1}$ and any integrable function $f$ defined on
$Q$, we let $|Q|$ denote the Lebesgue measure of $Q$ and let
$f_{Q}=\frac{1}{|Q|}\int_{Q}f(x)\d x$ denote the average of $f$ over
$Q$.  Analogous to the Morrey space $L^{p,\theta}(\Omega)$ and the
Campanato space ${\cal L}^{p,\theta}(\Omega)$, for $1\leq p<\infty$
we can define the parabolic Morrey space $L^{p,\theta}_{\rm
para}(\Omega_T)$ equipped with the norm
\begin{align*}
&\|f\|_{L^{p,\theta}_{\rm
para}(\Omega_T)}=\sup_{Q_R}\biggl(\frac{1}{R^{(n+2)\theta}}\int_{Q_R}|f(x)|^p\d
x\biggl)^\frac{1}{p}, \quad\quad\quad\quad 0\leq \theta\leq 1,
\end{align*}
and the parabolic Campanato space ${\cal L}^{p,\theta}_{\rm
para}(\Omega_T)$ of functions vanishing on the boundary
$\Gamma_T$, equipped with the norm
\begin{align*}
\|f\|_{{\cal
L}^{p,\theta}_{\rm para}(\Omega_T)}:
&=\sup_{Q_R(x_0,t_0)\cap\Omega_T}\biggl(\frac{1}{R^{(n+2)\theta}}\int_{Q_R(x_0,t_0)\cap\Omega_T}|f(x)|^p\d
x\biggl)^\frac{1}{p},\\
&~~~ +\sup_{Q_R(y_0,s_0)\cap\Omega_T}\biggl(\frac{1}{R^{(n+2)\theta}}\int_{Q_R(y_0,s_0)\cap\Omega_T}|f(x)-f_{Q_R}|^p\d
x\biggl)^\frac{1}{p},
\end{align*}
where the supremums above extend over all cylinders with
$x_0\in\partial\Omega$, $y_0\in\Omega$, $t_0,s_0\in(0,T]$ and
$0<R<R_\Omega$.

For any Banach space $X$ and time interval $(t_1,t_2)\subset\R$, we
denote by $L^p((t_1,t_2); X)$ the Bochner space equipped with the
norm
$$
\|f\|_{L^p((t_1,t_2);X)} =\left\{
\begin{array}{ll}
\displaystyle\biggl(\int_{t_1}^{t_2}\|f(t)\|_X^pdt\biggl)^\frac{1}{p},
& 1\leq p<\infty
,\\[15pt]
\displaystyle{\rm ess\,}\!\!\sup_{\!\!\!\!\!\!t\in
(t_1,t_2)}\|f(t)\|_X, & p=\infty.
\end{array}
\right.
$$

The importance of the (parabolic) Morrey spaces is that
$L^{p,\theta}(\Omega)$ translates just like
$L^{p/(1-\theta)}(\Omega)$, i.e. through the transformation $\tilde
f(y)=f(Ry)$ we have
\begin{align*}
&\|f\|_{L^{p,\theta}(B_R)} = R^{n(1-\theta)/p}\|\tilde
f\|_{L^{p,\theta}(B_1)},\\[-25pt]
\end{align*}
just like
\begin{align*}
&\|f\|_{L^{p/(1-\theta)}(B_R)} = R^{n(1-\theta)/p}\|\tilde
f\|_{L^{p/(1-\theta)}(B_1)},
\end{align*}
for any ball $B_R\subset\Omega$. Similarly, $L^{p,\theta}_{\rm
para}(\Omega_T)$ translates just like $L^{p/(1-\theta)}(\Omega_T)$. Therefore, $L^{p,\theta}(\Omega) $ and
$L^{p,\theta}_{\rm para}(\Omega_T)$ can be used as substitute for
$L^{p/(1-\theta)}(\Omega)$ and $L^{p/(1-\theta)}(\Omega_T)$, respectively, with lower order integrability. The
importance of the (parabolic) Campanato spaces includes: \\
(1)~ ${\cal L}^{p,1}(\Omega)$ are equivalent for all $1\leq
p<\infty$, i.e. ${\cal L}^{p,1}(\Omega)\cong \overline{\rm BMO}$;\\
(2)~ If $1<\theta<(n+p)/n$, then ${\cal L}^{p,\theta}(\Omega)\cong
C^{\alpha}_0(\overline\Omega)$ for $\alpha=n(\theta-1)/p$. \\
(3)~ If $1<\theta<(n+2+p)/(n+2)$, then ${\cal L}^{p,\theta}_{\rm
para}(\Omega_T)\cong C^{\alpha,\alpha/2}_0(\overline\Omega_T)$ for
$\alpha=(n+2)(\theta-1)/p$.

These properties of the Morrey and Campanato spaces can be found in
\cite{ChenYZParabEq,Troianiello}.

In this paper, we let $C_{p_1,p_2,\cdots,p_m}$ denote a generic
positive constant which depends on the parameters
$p_1,p_2,\cdots,p_m$.

\section{Main results}\label{dfjksl7879}
First, we establish a uniform-in-time BMO estimate and
a H\"{o}lder estimate for the solution of the
parabolic equation
\begin{align}
&\left\{
\begin{array}{ll} \displaystyle\frac{\partial u}{\partial
t}-\nabla\cdot(A\nabla u)=\nabla\cdot{\vec f} + f_0,
&\mbox{in}~~\Omega\times(0,T),\\[7pt]
u=g &\mbox{on} ~~\partial\Omega\times(0,T),\\[3pt]
u(x,0)=u_0(x)  &\mbox{for} ~~x\in\Omega ,
\end{array}
\right. \label{parabEq0}
\end{align}
where $\Omega$ is a bounded Lipschitz domain in $\R^n$ and
$A(x,t)=[A_{ij}(x,t)]_{n\times n}$ is a symmetric positive definite
measurable matrix function defined on $\R^{n+1}$ such that
\begin{align}\label{djfi}
K^{-1}|\xi|^2\leq \sum_{i,j=1}^n A_{ij}(x,t)\xi_i\xi_j \leq K|\xi|^2
,\quad\mbox{for all}~\xi\in\R^n
\end{align}
holds almost everywhere for $(x,t)\in\R^{n+1}$, where $K$ is a
positive constant.

\begin{theorem}\label{Mainresult}{\bf(BMO and H\"{o}lder estimates of parabolic equations)}~\\
{\it There exist positive constants $C$ and $\alpha_0\in (0,1)$
depending only on the elliptic constant $K$, the domain $\Omega$ and
the dimension $n$ $($independent of $T)$, such that the solution of
{\rm(\ref{parabEq0})} satisfies the BMO estimate
\begin{align}
&\|u\|_{L^\infty((0,T);\overline{\rm BMO})} \label{mainest1} \leq C(
\|f_0\|_{L^{1,n/(n+2)}(\Omega_T)}+\|{\vec
f}\|_{L^{2,n/(n+2)}(\Omega_T)}+\|u_0\|_{L^\infty(\Omega)}+
\|g\|_{L^\infty(\Gamma_T)}) .
\end{align}
If the compatibility condition $u_0(x)=g(x,0)$ for
$x\in\partial\Omega$ is satisfied, then we have
\begin{align}
&\|u\|_{C^{\alpha,\alpha/2}(\overline{\Omega}_T)}
\label{mainest2}\leq
C(\|f_0\|_{L^{1,(n+\alpha)/(n+2)}(\Omega_T)}+\|{\vec
f}\|_{L^{2,(n+2\alpha)/(n+2)}(\Omega_T)}
+\|u_0\|_{C^{\alpha}(\overline{\Omega})}+
\|g\|_{C^{\alpha,\alpha/2}(\Gamma_T)}) ,
\end{align}
for $0<\alpha\leq \alpha_0$.
}
\end{theorem}
The inequality (\ref{mainest1}) is new. A similar inequality as
(\ref{mainest2}) was proved in \cite{Yin}, where
$\|f_0\|_{L^{1,(n+\alpha)/(n+2)}(\Omega_T)}$ was replaced by
$\|f_0\|_{L^{2,(n-2+2\alpha)/(n+2)}(\Omega_T)}$. Note that
$L^{2,(n-2+2\alpha)/(n+2)}(\Omega_T)$ translates in the same way as
$L^{1,(n+\alpha)/(n+2)}(\Omega_T)$ under a scale transformation but
requires higher integrability. \bigskip

Secondly, by applying Theorem \ref{Mainresult}, we prove global
existence and uniqueness of a weak solution for the degenerate
thermistor problem
under the following physical hypotheses:\\[5pt]
{\bf(H1)} The thermal conductivity is a smooth function of
temperature and satisfies that
\begin{align*}
0<\inf_{s\geq r}\kappa(s)\leq \sup_{s\geq
r}\kappa(s)<\infty,\quad\mbox{for any fixed}~\, r>0 .
\end{align*}
{\bf(H2)} The electric resistivity $\rho(u)=1/\sigma(u)$
is a smooth function of temperature such that for some $p>0$ there holds
\begin{align}\label{dfsjljiowuo}
&C_{1,r}+C_{2,r}s^{p}\leq\rho(s)\leq C_{3,r}+C_{4,r}s^{p} \quad\mbox{$\forall~s\geq r>0$},
\end{align}
where $C_{i,r}$, $i=1,\cdots,5$, are some positive c
onstants (possibly depending on $r$).

Clearly, the hypotheses (H1)-(H2) are true for metals and some
semiconductors. In particular, the electric resistivity $\rho(u)$
can be any polynomials which are positive for $u> 0$. The hypotheses
(H1)-(H2) also imply that for any given $r>0$, $\sigma(s)$ is
bounded for $s\geq r$.

\begin{theorem}\label{mainthm}{\bf(Global well-posedness of the degenerate thermistor problem)}\\
{\it Let $\Omega$ be a bounded Lipschitz domain in $\R^n$ $(n=2,3)$ and let $q_0>n$. Assume that $u_0\in W^{1,q_0}(\Omega)$, $g\in L^\infty((0,T);W^{1,q_0}(\Omega))$, $\partial_tg\in L^\infty((0,T);L^{q_0}(\Omega))$,
$h\in L^\infty((0,T);W^{1,q_0}(\Omega))$,
with
$$\min_{(x,t)\in\Gamma_T}g(x,t)>0,\qquad
\min_{x\in
\Omega}\, u_0(x)
>0,
$$
and $g(x,0)=u_0(x)$ for $x\in\partial\Omega$. Then, under the hypothesis
{\rm (H1)-(H2)}, the initial-boundary value problem
{\rm (\ref{e-heat-1})-(\ref{BC})} admits a unique weak solution $(u,\phi)$ such that
\begin{align}
\begin{array}{ll}
u\in  C^{\alpha,\alpha/2}(\overline\Omega_T)\cap
L^p((0,T);W^{1,q}(\Omega)), \quad \phi\in  L^\infty((0,T);W^{1,q}(\Omega)) ,\\[3pt]
\partial_tu\in L^p((0,T);W^{-1,q}(\Omega)),
\end{array}
\label{sdfjhkl79031}
\end{align}
for some $q>n$, $0<\alpha<1$ and any $1<p<\infty$,  in the sense that the equations
\begin{align*}
&\int_0^T\int_{\Omega}\frac{\partial u}{\partial t}v\,\d x\d t
+\int_0^T\int_{\Omega}\kappa(u)\nabla
u\cdot \nabla v\,\d x\d t
=\int_0^T\int_{\Omega}\sigma(u)|\nabla\phi|^2v\,\d x\d t,\\
&\int_0^T\int_\Omega\sigma(u)\nabla\phi\cdot\nabla \varphi\,\d x\d t=0,\quad
\end{align*}
hold for any $v,\varphi\in L^2((0,T);H^1_0(\Omega))$.
}
\end{theorem}

Note that with the regularity (\ref{sdfjhkl79031}), the last equation above is equivalent to $$\int_\Omega\sigma(u)\nabla\phi\cdot\nabla \varphi\,\d x=0,\quad \forall~\varphi\in H^1_0(\Omega),\quad \mbox{a.e.}~t\in(0,T).$$


\section{BMO estimate of parabolic equations}\label{fjkesluio99}
\setcounter{equation}{0}

The solution of
(\ref{parabEq0}) can be decomposed into three parts, i.e.
the solution of  the following three problems:
\begin{align}
&\left\{
\begin{array}{ll}
\displaystyle\frac{\partial u}{\partial t}-\nabla\cdot(A\nabla u)=f_0,
&\mbox{in}~~\Omega\times(0,T),\\[7pt]
u=0 &\mbox{on} ~~\partial\Omega\times(0,T),\\[3pt]
u(x,0)=0 &\mbox{for} ~~x\in\Omega .
\end{array}
\right. \label{parabEq}\\[10pt]
&\left\{ \begin{array}{ll} \displaystyle\frac{\partial u}{\partial
t}-\nabla\cdot(A\nabla u)=\nabla\cdot{\vec f},
&\mbox{in}~~\Omega\times(0,T),\\[7pt]
u=0 &\mbox{on} ~~\partial\Omega\times(0,T),\\[3pt]
u(x,0)=0  &\mbox{for} ~~x\in\Omega ,
\end{array}
\right. \label{parabEq2}\\[10pt]
&\left\{ \begin{array}{ll} \displaystyle\frac{\partial u}{\partial
t}-\nabla\cdot(A\nabla u)=0,
&\mbox{in}~~\Omega\times(0,T),\\[7pt]
u=g &\mbox{on} ~~\partial\Omega\times(0,T),\\[3pt]
u(x,0)=u_0(x)  &\mbox{for} ~~x\in\Omega .
\end{array}
\right. \label{parabEq3}
\end{align}

From the maximum principle and the De Giorgi--Nash--Moser estimates, we know that there exist positive constants $C$ and $0<\alpha_0<1$ such that the solution of (\ref{parabEq3}) satisfies that
\begin{align*}
&\|u\|_{L^\infty(\Omega_T)} \leq \|g\|_{L^\infty(\Gamma_T)} +\|u_0\|_{L^\infty(\Omega)},\\
&\|u\|_{C^{\alpha,\alpha/2}(\overline{\Omega}_T)} \leq C
(\|g\|_{C^{\alpha,\alpha/2}(\overline\Gamma_T)}+\|u_0\|_{C^{\alpha}(\overline\Omega)}),
\end{align*}
for $0<\alpha<\alpha_0<1$ and $T>0$ (the second inequlaity above requires the compatability condition). To prove Theorem \ref{Mainresult}, it suffices to present estimates for the equations (\ref{parabEq})-(\ref{parabEq2}).

The rest part of this section is organized in the following way. In
Section \ref{secL11est}, we present local $L^1$ estimates for the
solution to (\ref{parabEq}). In Section \ref{glbest1}, we combine
the local $L^1$ estimates to derive a global BMO estimate based on
the equivalence of $\overline{\rm BMO}$ with the Campanato space
${\cal L}^{1,1}(\Omega)$. In Section \ref{sec4}, we establish the
BMO estimate for (\ref{parabEq2}) in terms of the Campanato space
${\cal L}^{2,1}(\Omega)$.

\subsection{Local $L^1$ estimates}\label{secL11est}
In this subsection, we present local $L^1$ estimates for the solution
of (\ref{parabEq}). The estimates obtained in this subsection will be used in Section \ref{glbest1} to derive a global BMO estimate uniformly with respect to time.

\begin{lemma}\label{Lemma001}
{\it Let $x_0\in\Omega$ and $0<t_0<T$. There exists
$\alpha_0\in(0,1)$ and $C>0$ such that if $u$ is the solution of
{\rm(\ref{parabEq})} in $Q_R=B_R(x_0)\times I_R$ with
$I_R=(t_0-R^2,t_0]$, then
\begin{align*}
&\max_{t\in I_\rho}\|u-u_{Q_\rho}\|_{L^1(B_\rho)} \leq
C\biggl(\frac{\rho}{R}\biggl)^{n+\alpha_0}\!\!\!\!\!\!\max_{t\in
I_R}\|u-\theta\|_{L^1(B_R)} +C\|f_0\|_{L^1(Q_R)}
\end{align*}
holds for any $0<\rho\leq R\leq \min({\rm
dist}(x_0,\partial\Omega),\sqrt{t_0})$ and any $\theta\in\R$, where
the constants $C$ and $\alpha_0$ depend only on $K$ and $n$. }
\end{lemma}
\noindent{\it Proof}~~~ First, we prove the lemma for
$\theta=0$. Let $\widetilde B_r=B_r(0)$, $\widetilde
I_r=(-r^2,0]$ and $\widetilde \Gamma_r=\partial\widetilde
B_r\times\widetilde I_r$. With any function $\xi$ defined on $Q_R$, we
associate a function $\tilde \xi(y,s)=\xi(x_0+Ry,t_0+R^2s)$ defined
on
 $\widetilde Q_1:=\widetilde B_1\times\widetilde I_1$. Then $\tilde u$ is a solution to the equation
$$
\frac{\partial \tilde u}{\partial s}-\nabla_y\cdot(\tilde A \nabla_y
\tilde u)=R^2\tilde f_0
$$
in $\widetilde Q_1$. Let $w$ be the solution of
\begin{align*}
\frac{\partial w}{\partial s}-\nabla_y\cdot(\tilde A \nabla_y
w)=R^2\tilde f_0
\end{align*}
with the boundary/initial condition $w=0$ on the parabolic
boundary $\partial_{\rm p}\widetilde Q_1$ and let $\bar w$ be the
solution of
\begin{align*}
\frac{\partial \bar w}{\partial s}-\nabla_y\cdot(\tilde A \nabla_y
\bar w)=R^2|\tilde f_0| 1_{\widetilde Q_1}
\end{align*}
in $\R^{n+1}$ with the initial condition $\bar w(y,0)\equiv 0$. By
the maximum principle, we know that
\begin{align*}
&|w(y,s)|\leq |\bar
w(y,s)|\\
&\leq \int_0^s\int_{\R^n}\frac{C}{(s-s')^{n/2}}e^{-\frac{|y-y'|^2}{C(s-s')}}R^2|\tilde
f_0(y',s')| 1_{\widetilde Q_1}(y',s')\d y'\d s'.
\end{align*}
Taking the $L^1(\widetilde B_1)$ norm with respect to $y$, we derive that
\begin{align*}
\|w\|_{L^\infty(\widetilde I_1;L^1(\widetilde B_1))}\leq
CR^2\|\tilde f_0\|_{L^1(\widetilde Q_1)} .
\end{align*}
We note that $v=\tilde u-\tilde u_{\widetilde Q_1}-w$ is the
solution of
\begin{align*}
&\frac{\partial v}{\partial s}-\nabla_y\cdot(\tilde  A \nabla_y v)=0
\end{align*}
in $\widetilde Q_1$, and by the De Giorgi--Nash estimates of parabolic
equations we know that there exists $\alpha_0\in(0,1)$ such that
for $\rho\in(0,1/2]$,
\begin{align*}
&\max_{t\in \widetilde I_\rho}\frac{1}{\rho^{n+\alpha_0}}\int_{\widetilde
B_\rho}|v-v_{\widetilde Q_\rho}|\d y \nn\\
&\leq C|v|_{C^{\alpha_0,\alpha_0/2}(\widetilde Q_{1/2})} \leq
C\|v\|_{L^1(\widetilde Q_1)} \leq C\max_{t\in\widetilde
I_1}\|v\|_{L^1(\widetilde B_1)}
 .
\end{align*}
Therefore,
\begin{align*}
&\max_{t\in\widetilde I_\rho}\|\tilde u-\tilde u_{\widetilde Q_\rho}\|_{L^1(\widetilde B_\rho)}\nn\\[3pt]
&\leq \max_{t\in\widetilde I_\rho}\|v-v_{\widetilde
Q_\rho}\|_{L^1(\widetilde B_\rho)}+\max_{t\in
\widetilde I_\rho}\|w-w_{\widetilde Q_\rho}\|_{L^1(\widetilde B_\rho)} \nn\\
&\leq C\rho^{n+\alpha_0}\max_{t\in\widetilde
I_1}\|v\|_{L^1(\widetilde B_1)}+C\max_{t\in\widetilde
I_1}\|w\|_{L^1(\widetilde B_1)}   \nn\\[3pt]
&\leq C\rho^{n+\alpha_0}\max_{t\in\widetilde I_1}\|\tilde u-\tilde
u_{\widetilde Q_1}\|_{L^1(\widetilde B_1)}+C\max_{t\in\widetilde
I_1}\|w\|_{L^1(\widetilde B_1)} \nn \\
&\leq C\rho^{n+\alpha_0}\max_{t\in\widetilde I_1}\|\tilde u-\tilde
u_{\widetilde Q_1}\|_{L^1(\widetilde B_1)}+CR^2\|\tilde
f_0\|_{L^1(\widetilde Q_1)}  \nn\\
&\leq C\rho^{n+\alpha_0}\max_{t\in\widetilde I_1}\|\tilde u\|_{L^1(\widetilde B_1)}+CR^2\|\tilde
f_0\|_{L^1(\widetilde Q_1)},
\end{align*}
where we have noted that
\begin{align*}
\|\widetilde u_{\widetilde Q_1}\|_{L^1(\widetilde B_1)}=\frac{|\widetilde B_1|}{|\widetilde Q_1|}\int_{\widetilde Q_1}|\widetilde u|\d x\d t\leq \max_{t\in \widetilde I_1}\|\widetilde u\|_{L^1(\widetilde B_1)} .
\end{align*}
Transforming back to the $(x,t)$-coordinates, we complete the proof of the Lemma.
for $\theta=0$.

Then we note that $u-\theta$ is also a
solution to the equation (\ref{parabEq}) in $Q_R$ for any $\theta\in\R$.\qed

Similarly, we can prove the following local $L^1$ estimates near the boundary $\partial_p\Omega_T$.

\begin{lemma}\label{Lemma002}
{\it Let $x_0\in\Omega$ and $t_0=0$. There exists $\alpha_0\in(0,1)$
and $C>0$ such that if $u$ is the solution of {\rm (\ref{parabEq})}
in $Q_R=B_R(x_0)\times\underline{I}_R$ with
$\underline{I}_R=[0,R^2]$, then
\begin{align*}
&\max_{t\in \underline{I}_\rho}\|u\|_{L^1(B_\rho)} \leq
C\biggl(\frac{\rho}{R}\biggl)^{n+\alpha_0}\!\!\!\!\!\!\max_{t\in
\underline{I}_R}\|u\|_{L^1(B_R)} +C\|f_0\|_{L^1(Q_R)}
\end{align*}
holds for any $0<\rho\leq R\leq \min({\rm
dist}(x_0,\partial\Omega),\sqrt{T})$, where the constants $C$ and
$\alpha_0$ depend only on $K$ and $n$. }
\end{lemma}

\begin{lemma}\label{Lemma003}
{\it Let $x_0\in\partial\Omega$ and $t_0>0$. There exists
$\alpha_0\in(0,1)$ and $C>0$ such that if $u$ is the solution of
{\rm (\ref{parabEq})} in $Q_R=B_R\times I_R$, with
$B_R=B_R(x_0)\cap\Omega$ and $I_R=(t_0-R^2,t_0]$, then
\begin{align*}
&\max_{t\in I_\rho}\|u\|_{L^1(B_\rho)}\leq
C\biggl(\frac{\rho}{R}\biggl)^{n+\alpha_0}\!\!\!\!\!\!\max_{t\in
I_R}\|u\|_{L^1(B_R)} +C\|f_0\|_{L^1(Q_R)}
\end{align*}
holds for any $0<\rho\leq R\leq \min(R_\Omega,\sqrt{t_0})$, where
the constants $C$ and $\alpha_0$ depend only on $K$, $n$ and
$\Omega$. }
\end{lemma}

\begin{lemma}\label{Lemma003}
{\it Let $x_0\in\partial\Omega$ and $t_0=0$. There exists
$\alpha_0\in(0,1)$ and $C>0$ such that if $u$ is the solution of
{\rm (\ref{parabEq})} in $Q_R=B_R\times I_R$, with
$B_R=B_R(x_0)\cap\Omega$ and $I_R=[0,R^2]$, then
\begin{align*}
&\max_{t\in I_\rho}\|u\|_{L^1(B_\rho)}\leq
C\biggl(\frac{\rho}{R}\biggl)^{n+\alpha_0}\!\!\!\!\!\!\max_{t\in
I_R}\|u\|_{L^1(B_R)} +C\|f_0\|_{L^1(Q_R)} ,
\end{align*}
holds for any $0<\rho\leq R\leq \min(R_\Omega,\sqrt{T})$, where the
constants $C$ and $\alpha_0$ depend only on $K$, $n$ and $\Omega$. }
\end{lemma}

The following simple lemma can be found in
\cite{ChenYZParabEq,Lieberman}, which is widely used for estimates
in terms of the Morrey and Campanato spaces.

\begin{lemma}\label{jkddlafhlae}
{\it Let $\varphi(\cdot)$ be a nonnegative and nondecreasing
function defined on $(0,R_0]$ and suppose that for any $0<\rho\leq R\leq
R_0$,
$$
\varphi(\rho)\leq
C_1\biggl(\frac{\rho}{R}\biggl)^{\gamma_1}\varphi(R)+C_2R^{\gamma_2}
,
$$
where $C_1$, $\gamma_1$ and $\gamma_2$ are nonnegative constants
such that $0<\gamma_2<\gamma_1$. Then
$$
\frac{1}{R^{\gamma_2}}\varphi(R)\leq
C_{\gamma_1,\gamma_2,C_1}\biggl(\frac{1}{R_0^{\gamma_2}}\varphi(R_0)+C_2\biggl)
.
$$
}
\end{lemma}\medskip

From the above lemmas, we obtain the following local $L^1$
estimates. \setcounter{proposition}{\value{lemma}}
\addtocounter{lemma}{1}
\begin{proposition}\label{dnflkqhklqewad}
{\it For $x_0\in\Omega$, $t_0>0$ and $Q_R=B_R(x_0)\times I_R$ with
$I_R=(t_0-R^2,t_0]$, we have
\begin{align*}
&\frac{1}{\rho^n}\|u-u_{Q_\rho}\|_{L^\infty((t_0-\rho^2,t_0);L^1(B_\rho))}\leq
C\biggl(\frac{1}{R^{n}}\|u\|_{L^\infty((t_0-R^2,t_0);L^1(B_{R}))}
+\|f_0\|_{L^{1,n/(n+2)}(\Omega_T)}\biggl) \nn
\end{align*}
for any $0<\rho\leq R\leq\min\big({\rm
dist}(x_0,\partial\Omega),\sqrt{t_0}\big)$. }
\end{proposition}

\addtocounter{lemma}{1}
\begin{proposition}\label{dnflkqhklqewad001}
{\it For $x_0\in\Omega$, $t_0=0$ and $Q_R=B_R(x_0)\times [0,R^2]$,
we have
\begin{align*}
&\frac{1}{\rho^n}\|u\|_{L^\infty((t_0-\rho^2,t_0);L^1(B_\rho))}\leq
C\biggl(\frac{1}{R^{n}}\|u\|_{L^\infty((t_0-R^2,t_0);L^1(B_{R}))}
+\|f_0\|_{L^{1,n/(n+2)}(\Omega_T)}\biggl)  \nn
\end{align*}
for any $0<\rho\leq R\leq\min\big({\rm
dist}(x_0,\partial\Omega),\sqrt{T}\big)$. }
\end{proposition}

\begin{proposition}\label{dnflkqhklqewad092}
{\it For $x_0\in\partial \Omega$, $t_0>0$ and
$Q_R=B_R(x_0)\cap\Omega\times I_R$ with $I_R=(t_0-R^2,t_0]$, we have
\begin{align*}
&\frac{1}{\rho^n}\|u\|_{L^\infty((t_0-\rho^2,t_0);L^1(B_\rho))}\leq
C\biggl(\frac{1}{R^{n}}\|u\|_{L^\infty((t_0-R^2,t_0);L^1(B_{R}))}
+\|f_0\|_{L^{1,n/(n+2)}(\Omega_T)}\biggl)  \nn
\end{align*}
for any $0<\rho\leq R\leq \min(R_\Omega,\sqrt{t_0})$.}
\end{proposition}

\begin{proposition}\label{dnflkqhklqewad03}
{\it For $x_0\in\partial \Omega$, $t_0=0$ and
$Q_R=B_R(x_0)\cap\Omega\times [0,R^2]$, we have
\begin{align*}
&\frac{1}{\rho^n}\|u\|_{L^\infty((0,\rho^2);L^1(B_\rho))}\leq
C\biggl(\frac{1}{R^{n}}\|u\|_{L^\infty((0,R^2);L^1(B_{R}))}
+\|f_0\|_{L^{1,n/(n+2)}(\Omega_T)}\biggl) \nn
\end{align*}
for any $0<\rho\leq R\leq \min(R_\Omega,\sqrt{T})$.}
\end{proposition}

\subsection{BMO estimates via ${\cal L}^{1,1}$ }\label{glbest1}
We combine the local $L^1$ estimates obtained in the last subsection
to derive a global BMO estimate of $u$, uniform with respect to time.

\begin{proposition}\label{njdfklahkeas}
{\it The Propositions
{\rm\ref{dnflkqhklqewad}--\ref{dnflkqhklqewad03}} imply that the solution of {\rm (\ref{parabEq})} satisfies that
\begin{align}\label{BMOest2df}
\|u\|_{L^\infty((0,T);{\rm \overline{BMO}})} \leq
C\|f_0\|_{L^{1,n/(n+2)}(\Omega_T)},
\end{align}
where $C$ depends only on $K$, $n$ and $\Omega$
(independent of $T$). }
\end{proposition}
\noindent{\it Proof}~~~ Set $M=\|f_0\|_{L^{1,n/(n+2)}(\Omega_T)}$.

First, we prove the proposition for $T\geq
R_\Omega^2$.
We shall prove that for $R< R_\Omega/2$ and any set
$B_R=B_R(x_0)\cap\Omega$ with some point $x_0\in\overline\Omega$ and
$\delta={\rm dist}(x_0,\partial\Omega)$, the following estimates
hold:
\begin{align}\label{estkeaklael1}
\left\{
\begin{array}{ll}
\frac{1}{R^n}\|u\|_{L^\infty((0,T);L^1(B_R))}\leq
C\big(\|u\|_{L^\infty((0,T);L^1(\Omega))} +M\big) , &
\mbox{if}~~\delta\leq  R,\\[5pt]
\frac{1}{R^n}\|u-u_{B_R}\|_{L^\infty((0,T);L^1(B_R))}\leq
C\big(\|u\|_{L^\infty((0,T);L^1(\Omega))} +M\big) , &
\mbox{if}~~\delta> R,
\end{array}
\right.
\end{align}
%
{\it Case} 1: $\delta \leq R$.~~ In this case, there exists a region
$B_{2R}=B_{2R}(y_0)\cap\Omega$ with some $y_0\in\partial\Omega$ such
that $B_R\subset B_{2R}$ and so, for any given $t_0\in[0,T]$,
\begin{align}\label{fjekqhlqt}
&\|u(\cdot,t_0)\|_{L^1(B_R)}\leq \|u(\cdot,t_0)\|_{L^1(B_{2R})} .
\end{align}
Now if $t_0\leq 4R^2$, then by Proposition \ref{dnflkqhklqewad03},
\begin{align}
&\frac{1}{R^n}\|u(\cdot,t_0)\|_{L^1(B_{2R})}\leq \|u\|_{L^\infty((0,4R^2);L^1(B_{2R}))}\nn\\
&\leq
C\biggl(\frac{1}{R_\Omega^{n}}\|u\|_{L^\infty((0,R^2_\Omega);L^1(B_{R_\Omega}))}
+M\biggl) .
\end{align}
Otherwise, $t_0>4R^2$ and by Proposition \ref{dnflkqhklqewad092}, for
$R_0=\min(\sqrt{t_0},R_\Omega)$ and $R_m=\max(\sqrt{t_0},R_\Omega)$ we have
\begin{align*}
&\frac{1}{R^n}\|u\|_{L^\infty((t_0-4R^2,t_0);L^1(B_{2R}))} \\
&\leq
C\biggl(\frac{1}{R_0^{n}}\|u\|_{L^\infty((t_0-R_0^2,t_0);L^1(B_{R_0}))}
+M\biggl) \\
&= \left\{
\begin{array}{ll}
C\biggl(\frac{1}{t_0^{n/2}}\|u\|_{L^\infty((0,t_0);L^1(B_{\sqrt{t_0}}))}
+M\biggl), & \sqrt{t_0}<R_\Omega,\\[5pt]
C\biggl(\frac{1}{R_\Omega^{n}}\|u\|_{L^\infty((t_0-R_\Omega^2,t_0);L^1(B_{R_\Omega}))}
+M\biggl), & \sqrt{t_0}\geq R_\Omega .
\end{array}
\right.\\
&\leq \left\{
\begin{array}{ll}
C\biggl(\frac{1}{R_\Omega^{n}}\|u\|_{L^\infty((0,R_\Omega);L^1(B_{R_\Omega}))}
+M\biggl) , & \sqrt{t_0}<R_\Omega, \quad\mbox{(by Proposition \ref{dnflkqhklqewad03})} \\[5pt]
C\biggl(\frac{1}{R_\Omega^{n}}\|u\|_{L^\infty((t_0-R_\Omega^2,t_0);L^1(B_{R_\Omega}))}
+M\biggl), & \sqrt{t_0}\geq R_\Omega .
\end{array}
\right.
\end{align*}
To conclude, for $\delta\leq R$ and $t_0\in[0,T]$ we have
\begin{align}\label{fdjkaletehw}
\frac{1}{R^n}\|u(\cdot,t_0)\|_{L^1(B_R)}\leq
\frac{C}{R_\Omega^{n}}\|u\|_{L^\infty((0,T);L^1(\Omega))}+CM.
\end{align}

\noindent{\it Case} 2: $\delta>R$. ~ In this case, we set
$R_0=\min(\delta,\sqrt{t_0},R_\Omega)$. Then Proposition
\ref{dnflkqhklqewad} implies that
\begin{align*}
&\frac{1}{R^n}\|u-u_{B_R}\|_{L^\infty((t_0-R^2,t_0);L^1(B_R))} \leq \frac{1}{R^n}\|u-u_{Q_R}\|_{L^\infty((t_0-R^2,t_0);L^1(B_R))} \\
&\leq
C\biggl(\frac{1}{R_0^{n}}\|u\|_{L^\infty((t_0-R_0^2,t_0);L^1(B_{R_0}))}
+M\biggl)\\
&= \left\{
\begin{array}{ll}
C\biggl(\frac{1}{R_\Omega^{n}}\|u\|_{L^\infty((t_0-R_\Omega^2,t_0);L^1(B_{R_\Omega}))}
+M\biggl) , & \mbox{if}~~R_\Omega\leq\min(\delta,\sqrt{t_0})\\
C\biggl(\frac{1}{\delta^{n}}\|u\|_{L^\infty((t_0-\delta^2,t_0);L^1(B_\delta))}
+M\biggl) , & \mbox{else if}~~\delta\leq\min(\sqrt{t_0},R_\Omega)\\
C\biggl(\frac{1}{R_0^{n}}\|u\|_{L^\infty((0,R_0^2);L^1(B_{R_0}))}
+M\biggl) , & \mbox{else }~~t_0=R_0^2
\end{array}
\right.\\
&\leq \left\{
\begin{array}{ll}
\frac{C}{R_\Omega^{n}}\|u\|_{L^\infty((0,T);L^1(\Omega))} +CM,
~~~~~~~~~~~~~~~~~~\mbox{if}~~R_\Omega\leq\min(\delta,\sqrt{t_0})\\[8pt]
\frac{C}{R_\Omega^{n}}\|u\|_{L^\infty((0,T);L^1(\Omega))}  +CM,
~~~~~~~~~~~~~~~~~~\mbox{else if}~~\delta\leq\min(\sqrt{t_0},R_\Omega) ~~~\mbox{by}~(\ref{fdjkaletehw})\\[8pt]
\left\{
\begin{array}{ll}
C\biggl(\frac{1}{\delta^{n}}\|u\|_{L^\infty((0,\delta^2);L^1(B_\delta))}
+M\biggl), & \mbox{else if}~~\delta\leq R_\Omega\\[5pt]
C\biggl(\frac{1}{R_\Omega^{n}}\|u\|_{L^\infty((0,R_\Omega^2);L^1(B_{R_\Omega}))}
+M\biggl), & \mbox{else if}~~\delta>R_\Omega
\end{array}
\right.
\quad\mbox{(by Proposition \ref{dnflkqhklqewad03})}
\end{array}
\right.
\\[5pt]
&\leq \frac{C}{R_\Omega^{n}}\|u\|_{L^\infty((0,T);L^1(\Omega))}
+CM,\quad \mbox{again by (\ref{fdjkaletehw})} .
\end{align*}

So far we have proved (\ref{estkeaklael1}). Once we note that
$\|u\|_{L^\infty((0,T);L^1(\Omega))}\leq
C\|f_0\|_{L^1(\Omega_T)}$, we derive (\ref{BMOest2df}) from (\ref{estkeaklael1}).

Secondly, we prove the proposition for $0<T<R_\Omega$. In this case, we
consider the solution $\hat{u}$ of the equation
\begin{align}
\frac{\partial\hat u}{\partial t}-\nabla\cdot( A \nabla\hat u )=\hat
f_0
\end{align}
in the domain $\Omega_{R_\Omega}=\Omega\times(0,R_\Omega)$ with the boundary and initial conditions $\hat
u =0$ on $\partial\Omega\times(0,R_\Omega)$ and $\hat u (x,0)=0$ for
$x\in\Omega$, where
$$
\hat f_0(x,t)=\left\{
\begin{array}{ll}
f_0(x,t), &\mbox{for}~~t\in(0,T),\\[5pt]
0,      &\mbox{for}~~t\in(T,R_\Omega) .
\end{array}\right.
$$
Check that
\begin{align*}
&\|\hat f_0\|_{L^{1,n/(n+2)}(\Omega_{R_\Omega})}\leq
C\|f_0\|_{L^{1,n/(n+2)}(\Omega_T)},\\
&\|\hat f_0\|_{L^1(\Omega_{R_\Omega})}\leq C\|f_0\|_{L^1(\Omega_T)} ,\\
&\|u\|_{L^\infty((0,T);{\rm \overline{BMO}})}\leq \|\hat
u\|_{L^\infty((0,R_\Omega);{\rm \overline{BMO}})},
\end{align*}
where the constant $C$ does not depend on $T$ (as $T\rightarrow 0$). Then we
apply the inequality (\ref{BMOest2df}) to $\hat u$ with $T=R_\Omega$.\qed

\subsection{BMO estimates via ${\cal L}^{2,1}$}\label{sec4}
\setcounter{equation}{0}

In this section, we present estimates for the solution of
(\ref{parabEq2}). The idea is similar as Section \ref{glbest1}. From the proof of the following lemma we can see the main
difference between the current subsection and the last subsection.

\begin{lemma}\label{Lemma004}
{\it Let $x_0\in\Omega$ and $0<t_0<T$. There exists
$\alpha_0\in(0,1)$ and $C>0$ such that if $u$ is the solution to
{\rm (\ref{parabEq2})} in $Q_R=B_R(x_0)\times I_R$ with
$I_R=(t_0-R^2,t_0]$, then
\begin{align*}
&\max_{t\in I_\rho}\|u-u_{Q_\rho}\|_{L^2(B_\rho)}^2 \leq
C\biggl(\frac{\rho}{R}\biggl)^{n+2\alpha_0}
\!\!\!\!\!\!\max_{t\in
I_R}\|u-\theta\|_{L^2(B_R)}^2 +C\|\vec{f}\|_{L^2(Q_R)}^2
\end{align*}
holds for any $0<\rho\leq R\leq \min({\rm
dist}(x_0,\partial\Omega),\sqrt{t_0})$ and any $\theta\in\R$, where
$C$ depends only on $K$ and $n$. }
\end{lemma}
\noindent{\it Proof}~~~ Let $\widetilde B_r=B_r(0)$, $\widetilde
I_r=(-r^2,0]$ and $\widetilde \Gamma_r=\partial\widetilde
B_r\times\widetilde I_r$. With any function $w$ defined on $Q_R$, we
associate a function $\tilde \xi(y,s)=\xi(x_0+Ry,t_0+R^2s)$ defined
on
 $\widetilde Q_1:=\widetilde B_1\times\widetilde I_1$. Then $\tilde u$ is a solution of the equation
$$
\frac{\partial \tilde u}{\partial s}-\nabla_y\cdot(\tilde A \nabla_y
\tilde u)=R\nabla_y\cdot \tilde f
$$
in $\widetilde Q_1$. Let $w$ be the solution of
\begin{align*}
\frac{\partial w}{\partial s}-\nabla_y\cdot(\tilde A \nabla_y
w)=R\nabla_y\cdot \tilde f
\end{align*}
with the initial and boundary condition $w=0$ on the parabolic
boundary $\partial_{\rm p}\widetilde Q_1$. Multiplying the above
equation by $w$ and integrating the result over $\widetilde Q_1$, we
obtain that
\begin{align*}
&\|w\|_{L^\infty(\widetilde I_1;L^2(\widetilde B_1))}\leq
CR\|\tilde f\|_{L^2(\widetilde Q_1)}
\end{align*}
On the other hand, we observe that $v=\tilde u-\tilde u_{\widetilde
Q_1}-w$ is the solution of
\begin{align*}
&\frac{\partial v}{\partial s}-\nabla_y\cdot(\tilde  A \nabla_y v)=0
\end{align*}
in $\widetilde Q_1$. By the De Giorgi--Nash estimates of parabolic
equations, we know that there exists $\alpha_0\in(0,1)$ such that
for $\rho\in(0,1/2]$,
\begin{align*}
&\max_{t\in I_\rho}\frac{1}{\rho^{n+2\alpha_0}}\int_{\widetilde
B_\rho}|v-v_{\widetilde Q_\rho}|^2\d y \nn\\
&\leq C|v|_{C^{\alpha_0,\alpha_0/2}(\widetilde Q_{1/2})}^2 \leq
C\|v\|_{L^2(\widetilde Q_1)}^2 \leq C\max_{t\in\widetilde
I_1}\|v\|_{L^2(\widetilde B_1)}^2
 .
\end{align*}
Therefore,
\begin{align*}
&\max_{t\in\widetilde I_\rho}\|\tilde u-\tilde u_{\widetilde Q_\rho}\|_{L^2(\widetilde B_\rho)}^2\nn\\[3pt]
&\leq C\max_{t\in\widetilde I_\rho}\|v-v_{\widetilde
Q_\rho}\|_{L^2(\widetilde B_\rho)}^2+C\max_{t\in
\widetilde I_\rho}\|w-w_{\widetilde Q_\rho}\|_{L^2(\widetilde B_\rho)}^2 \nn\\
&\leq C\rho^{n+2\alpha_0}\max_{t\in\widetilde
I_1}\|v\|_{L^2(\widetilde B_1)}^2+C\max_{t\in\widetilde
I_1}\|w\|_{L^2(\widetilde B_1)}^2   \nn\\[3pt]
&\leq C\rho^{n+2\alpha_0}\max_{t\in\widetilde I_1}\|\tilde u-\tilde
u_{\widetilde Q_1}\|_{L^2(\widetilde B_1)}^2+C\max_{t\in\widetilde
I_1}\|w\|_{L^2(\widetilde B_1)}^2 \nn \\
&\leq C\rho^{n+2\alpha_0}\max_{t\in\widetilde I_1}\|\tilde u-\tilde
u_{\widetilde Q_1}\|_{L^2(\widetilde B_1)}^2+CR^2\|\tilde f\|_{L^2(\widetilde Q_1)}^2  \nn\\
&\leq C\rho^{n+2\alpha_0}\max_{t\in\widetilde I_1}
\|\tilde u\|_{L^2(\widetilde B_1)}^2+CR^2\|\tilde f\|_{L^2(\widetilde Q_1)}^2  .
\end{align*}
Transforming back to the $(x,t)$-coordinates, we complete the proof
of the
Lemma for $\theta=0$.

Then we note that $u-\theta$ is also a
solution to the equation (\ref{parabEq2}) in $Q_R$ for any $\theta\in\R$.\qed\medskip

In a similar way, we can prove the following lemmas and propositions.

\begin{lemma}\label{Lemma005}
{\it Let $x_0\in\Omega$ and $t_0=0$. There exists $\alpha_0\in(0,1)$
and $C>0$ such that if $u$ is the solution of {\rm (\ref{parabEq2})}
in $Q_R=B_R(x_0)\times \underline{I}_R$ with
$\underline{I}_R=[0,R^2]$, then
\begin{align*}
&\max_{t\in \underline{I}_\rho}\|u\|_{L^2(B_\rho)}^2 \leq
C\biggl(\frac{\rho}{R}\biggl)^{n+2\alpha_0}\!\!\!\!\!\!\max_{t\in
\underline{I}_R}\|u\|_{L^2(B_R)}^2  +C\|\vec{f}\|_{L^2(Q_R)}^2
\end{align*}
holds for any $0<\rho\leq R\leq \min({\rm
dist}(x_0,\partial\Omega),\sqrt{T})$,
where $C$ and $\alpha_0$ depend only on $K$ and
$n$. }
\end{lemma}

\begin{lemma}\label{Lemma006}
{\it Let $x_0\in\partial\Omega$ and $t_0>0$. There exists
$\alpha_0\in(0,1)$ and $C>0$ such that if $u$ is the solution of {\rm
(\ref{parabEq2})} in $Q_R=B_R\times I_R$ with
$B_R=B_R(x_0)\cap\Omega$ and $I_R=(t_0-R^2,t_0]$, then
\begin{align*}
&\max_{t\in I_\rho}\|u\|_{L^2(B_\rho)}^2\leq
C\biggl(\frac{\rho}{R}\biggl)^{n+2\alpha_0}\max_{t\in
I_R}\|u\|_{L^2(B_R)}^2 +C\|\vec{f}\|_{L^2(Q_R)}^2
\end{align*}
holds for any  $0<\rho\leq R\leq
\min(R_\Omega,\sqrt{t_0})$, where $C$ and $\alpha_0$ depend only on $K$, $n$ and $\Omega$. }
\end{lemma}

\begin{lemma}\label{Lemma007}
{\it Let $x_0\in\partial\Omega$ and $t_0=0$. There exists
$\alpha_0\in(0,1)$ and $C>0$ such that if $u$ is the solution to
{\rm(\ref{parabEq2})} in $Q_R=B_R\times\underline{ I}_R$, with
$B_R=B_R(x_0)\cap\Omega$ and $\underline{I}_R=[0,R^2]$, then
\begin{align*}
&\max_{t\in \underline{I}_\rho}\|u\|_{L^2(B_\rho)}^2\leq
C\biggl(\frac{\rho}{R}\biggl)^{n+2\alpha_0}\!\!\!\!\!\!\max_{t\in
\underline{I}_R}\|u\|_{L^2(B_R)}^2 +C\|\vec{f}\|_{L^2(Q_R)}^2
\end{align*}
holds for any $0<\rho\leq R\leq
\min(R_\Omega,\sqrt{T})$, where $C$ and $\alpha_0$ depend only on $K$, $n$ and $\Omega$. }
\end{lemma}

From the above lemmas, using Lemma \ref{jkddlafhlae} we can derive
the following results concerning the solution of (\ref{parabEq2}).

\setcounter{proposition}{\value{lemma}} \addtocounter{lemma}{1}
\begin{proposition}\label{dnflkqhklqewad2}
{\it For $x_0\in\Omega$, $t_0>0$, $Q_R=B_R(x_0)\times(t_0-R^2,t_0]$
and $0<\rho\leq R\leq\min\big({\rm
dist}(x_0,\partial\Omega),\sqrt{t_0}\big)$, we have
\begin{align*}
&\|u-u_{Q_\rho}\|_{L^\infty((t_0-\rho^2,t_0);L^2(B_\rho))}^2\leq
C\biggl(\frac{1}{R^{n}}\|u\|_{L^\infty((t_0-R^2,t_0);L^2(B_{R}))}^2
+\|\vec{f}\|_{L^{2,n/(n+2)}(\Omega_T)}^2\biggl)\rho^n .\nn
\end{align*}
}
\end{proposition}


\addtocounter{lemma}{1}
\begin{proposition}\label{dnflkqhklqewad12}
{\it For $x_0\in\Omega$, $t_0=0$, $Q_R=B_R(x_0)\times[0,R^2]$ and
$0<\rho\leq R\leq\min\big({\rm
dist}(x_0,\partial\Omega),\sqrt{T}\big)$, we have
\begin{align*}
&\|u\|_{L^\infty((0,\rho^2);L^2(B_\rho))}^2\leq
C\biggl(\frac{1}{R^{n}}\|u\|_{L^\infty((0,R^2);L^2(B_{R}))}^2
+\|\vec{f}\|_{L^{2,n/(n+2)}(\Omega_T)}^2\biggl)\rho^n . \nn
\end{align*}
}
\end{proposition}

\begin{proposition}\label{dnflkqhklqewad22}
{\it For $x_0\in\partial \Omega$, $t_0>0$,
$Q_R=B_R(x_0)\cap\Omega\times(t_0-R^2,t_0]$ and
$0<\rho<R\leq\min(R_\Omega,\sqrt{t_0})$, we have
\begin{align*}
&\|u\|_{L^\infty((t_0-\rho^2,t_0);L^2(B_\rho))}^2\leq
C\biggl(\frac{1}{R^{n}}\|u\|_{L^\infty((t_0-R^2,t_0);L^2(B_{R}))}^2
+\|\vec{f}\|_{L^{2,n/(n+2)}(\Omega_T)}^2\biggl)\rho^n . \nn
\end{align*}
}
\end{proposition}

\begin{proposition}\label{dnflkqhklqewad32}
{\it For $x_0\in\partial \Omega$, $t_0=0$,
$Q_R=B_R(x_0)\cap\Omega\times[0,R^2]$ and $0<\rho<R\leq
\min(R_\Omega,\sqrt{T})$, we have
\begin{align*}
&\|u\|_{L^\infty((0,\rho^2);L^2(B_\rho))}^2\leq C\biggl(\frac{1}{R^{n}}\|u\|_{L^\infty((0,R^2);L^2(B_{R}))}^2
+\|\vec{f}\|_{L^{2,n/(n+2)}(\Omega_T)}^2\biggl)\rho^n . \nn
\end{align*}
}
\end{proposition}

With the above propositions and following the outline of Section
\ref{glbest1}, we can prove the global BMO estimate below.
\begin{proposition}\label{njdfklahkeas2222}
{\it The Propositions
{\rm\ref{dnflkqhklqewad2}--\ref{dnflkqhklqewad32}} imply that the solution of {\rm (\ref{parabEq2})} satisfies that
\begin{align}\label{BMOest22}
\|u\|_{L^\infty((0,T);{\rm \overline{BMO}})} \leq C\|\vec{f}\|_{L^{2,n/(n+2)}(\Omega_T)},
\end{align}
where $C$ depends only on $K$, $n$ and $\Omega$
$($independent of $T)$. }
\end{proposition}

\section{H\"{o}lder estimate of parabolic equations}\label{sec5}
\setcounter{equation}{0}

In this section, we list the propositions to be used in deriving (\ref{mainest2}). We omit the proof of these propositions, since it is very similar as the last section,
The reason we keep these propositions in this section is that some of them are also used in the next section to prove global well-posedness of the degenerate thermistor problem.

There exist positive constants $\alpha_0$ and $C$ such that the following propositions hold.

\setcounter{proposition}{\value{lemma}} \addtocounter{lemma}{1}
\begin{proposition}\label{dnflkqhklqewad3}
{\it For $x_0\in\Omega$, $t_0>0$, $Q_R=B_R(x_0)\times(t_0-R^2,t_0]$,
$0<2\rho\leq R\leq\min\big({\rm
dist}(x_0,\partial\Omega),\sqrt{t_0}\big)$ and $0<\alpha<\alpha_0$,
the solution of {\rm (\ref{parabEq})} satisfies that
\begin{align*}
&\frac{1}{\rho^{n+2+\alpha}}\|u-u_{Q_\rho}\|_{L^1(Q_\rho)}\leq
C\biggl(\frac{1}{R^{n+2+\alpha}}\|u-\theta\|_{L^1(Q_{R})}
+\frac{1}{R^{n+\alpha}}\|f_0\|_{L^1(Q_R)}\biggl), \nn
\end{align*}
where $\theta$ is an arbitrary constant. }
\end{proposition}

\addtocounter{lemma}{1}
\begin{proposition}\label{dnflkqhklqewad13}
{\it For $x_0\in\Omega$, $t_0=0$, $Q_R=B_R(x_0)\times[0,R^2]$,
$0<\rho\leq R\leq\min\big({\rm
dist}(x_0,\partial\Omega),\sqrt{T}\big)$ and $0<\alpha<\alpha_0$,
the solution of {\rm (\ref{parabEq})} satisfies that
\begin{align*}
&\frac{1}{\rho^{n+2+\alpha}}\|u\|_{L^1(Q_\rho)}\leq C\biggl(\frac{1}{R^{n+2+\alpha}}\|u\|_{L^1(Q_{R})}
+\frac{1}{R^{n+\alpha}}\|f_0\|_{L^1(Q_R)}\biggl)
. \nn
\end{align*}
}
\end{proposition}

\begin{proposition}\label{dnflkqhklqewad23}
{\it For $x_0\in\partial \Omega$, $t_0>0$,
$Q_R=B_R(x_0)\cap\Omega\times(t_0-R^2,t_0]$, $0<\rho\leq
R\leq\min(R_\Omega,\sqrt{t_0})$ and $0<\alpha<\alpha_0$, the
solution of {\rm (\ref{parabEq})} satisfies that
\begin{align*}
&\frac{1}{\rho^{n+2+\alpha}}\|u\|_{L^1(Q_\rho)}\leq
C\biggl(\frac{1}{R^{n+2+\alpha}}\|u\|_{L^1(Q_{R})}
+\frac{1}{R^{n+\alpha}}\|f_0\|_{L^1(Q_R)}\biggl) . \nn
\end{align*}
}
\end{proposition}

\begin{proposition}\label{dnflkqhklqewad33}
{\it For $x_0\in\partial \Omega$, $t_0=0$,
$Q_R=B_R(x_0)\cap\Omega\times[0,R^2]$, $0<\rho\leq R\leq
\min(R_\Omega,\sqrt{T})$  and $0<\alpha<\alpha_0$, the solution of
{\rm (\ref{parabEq})} satisfies that
\begin{align*}
&\frac{1}{\rho^{n+2+\alpha}}\|u\|_{L^1(Q_\rho)}\leq
C\biggl(\frac{1}{R^{n+2+\alpha}}\|u\|_{L^1(Q_{R})}
+\frac{1}{R^{n+\alpha}}\|f_0\|_{L^1(Q_R)}\biggl)
. \nn
\end{align*}
}
\end{proposition}

With the above propositions and following the outline of Section
\ref{glbest1}, we can derive the following estimate in terms of the Campanato space.
\begin{proposition}\label{njdfklahkeas3333}
{\it The solution of {\rm(\ref{parabEq})} satisfies that
\begin{align}\label{BMOest333}
\|u\|_{{\cal L}^{1,1+\alpha/(n+2)}_{\rm para}(\Omega_T)} \leq
 C\|f_0\|_{L^{1,(n+\alpha)/(n+2)}(\Omega_T)},
\end{align}
where $C$ depends only on $K$, $n$ and $\Omega$
$($independent of $T)$. }
\end{proposition}

The local and global estimates in ${\cal L}^{2,\theta}_{\rm para}(\Omega_T)$
follow in a similar way. To conclude, we have

\begin{proposition}\label{dnflkqaghu}
{\it For $x_0\in\Omega$, $t_0>0$, $Q_R=B_R(x_0)\times(t_0-R^2,t_0]$,
$0<2\rho\leq R\leq\min\big({\rm
dist}(x_0,\partial\Omega),\sqrt{t_0}\big)$ and $0<\alpha<\alpha_0$,
the solution of {\rm (\ref{parabEq2})} satisfies that
\begin{align*}
&\frac{1}{\rho^{n+2+2\alpha}}
\|u-u_{Q_\rho}\|_{L^2(Q_\rho)}^2\leq
C\biggl(\frac{1}{R^{n+2+2\alpha}}
\|u-\theta\|_{L^2(Q_{R})}^2
+\frac{1}{R^{n+2\alpha}}\|\vec{f}\|_{L^2(Q_R)}^2
\biggl), \nn
\end{align*}
where $\theta$ is an arbitrary constant. }
\end{proposition}

\addtocounter{lemma}{1}
\begin{proposition}\label{dnflkqas13}
{\it For $x_0\in\Omega$, $t_0=0$,
$Q_R=B_R(x_0)\cap\Omega\times[0,R^2]$, $0<\rho\leq
R\leq\min\big({\rm dist}(x_0,\partial\Omega),\sqrt{T}\big)$ and
$0<\alpha<\alpha_0$, the solution of {\rm (\ref{parabEq2})}
satisfies that
\begin{align*}
&\frac{1}{\rho^{n+2+2\alpha}}
\|u\|_{L^2(Q_\rho)}^2\leq C\biggl(\frac{1}{R^{n+2+2\alpha}}
\|u\|_{L^2(Q_{R})}^2
+\frac{1}{R^{n+2\alpha}}\|\vec{f}\|_{L^2(Q_R)}^2
\biggl)
. \nn
\end{align*}
}
\end{proposition}

\begin{proposition}\label{dnflkqhksgd23}
{\it For $x_0\in\partial \Omega$, $t_0>0$,
$Q_R=B_R(x_0)\cap\Omega\times(t_0-R^2,t_0]$,
$0<\rho<R\leq\min(R_\Omega,\sqrt{t_0})$ and $0<\alpha<\alpha_0$, the
solution of {\rm (\ref{parabEq2})} satisfies that
\begin{align*}
&\frac{1}{\rho^{n+2+2\alpha}}\|u\|_{L^1(Q_\rho)}\leq
C\biggl(\frac{1}{R^{n+2+2\alpha}}\|u\|_{L^2(Q_{R})}^2
+\frac{1}{R^{n+2\alpha}}\|\vec{f}\|_{L^2(Q_R)}^2
\biggl) . \nn
\end{align*}
}
\end{proposition}

\begin{proposition}\label{dnflkqsady63}
{\it For $x_0\in\partial \Omega$, $t_0=0$,
$Q_R=B_R(x_0)\cap\Omega\times[0,R^2]$, $0<\rho\leq R\leq
\min(R_\Omega,\sqrt{T})$ and $0<\alpha<\alpha_0$, the solution of
{\rm (\ref{parabEq2})} satisfies that
\begin{align*}
&\frac{1}{\rho^{n+2+2\alpha}}\|u\|_{L^2(Q_\rho)}^2\leq
C\biggl(\frac{1}{R^{n+2+2\alpha}}\|u\|_{L^2(Q_{R})}^2
+\frac{1}{R^{n+2\alpha}}\|\vec{f}\|_{L^2(Q_R)}^2
\biggl)
. \nn
\end{align*}
}
\end{proposition}

\begin{proposition}\label{njdfklah5yj}
{\it The solution of {\rm(\ref{parabEq2})} satisfies that
\begin{align}\label{BMOest333}
\|u\|_{{\cal L}^{2,1+2\alpha/(n+2)}_{\rm para}(\Omega_T)} \leq
C\|\vec{f}\|_{L^{2,(n+2\alpha)/(n+2)}(\Omega_T)}
,
\end{align}
where $C$ depends only on $K$, $n$ and $\Omega$
$($independent of $T)$. }
\end{proposition}

Proposition \ref{njdfklahkeas3333} and Proposition
\ref{njdfklah5yj} imply the global H\"{o}lder estimate
(\ref{mainest2}).

\bigskip

\section{The degenerate thermistor problem}\label{fsdjkljl78687}
\setcounter{equation}{0}

In this section, we prove Theorem \ref{mainthm} concerning
global well-posedness of the degenerate thermistor problem.
Before we prove the theorem, we introduce some lemmas to be used.

\subsection{Preliminaries}

\begin{lemma}\label{BMOandrecipr}
{\it Let $p>0$. If $u\in$ BMO$(\R^n)$, $u\geq 0$, and $C_1+C_2|s|^p\leq\rho(s)\leq C_3+C_4|s|^p$ for $s\geq 0$, then $\rho(u)$ is a $A_2$ weight in the sense that
\begin{align*}
\biggl(\frac{1}{|B|}\int_{B}\rho(u)\d
x\biggl)\biggl(\frac{1}{|B|}\int_{B}\frac{1}{\rho(u)}\d x\biggl)\leq
C
\end{align*}
for any ball $B\subset\R^n$, where the constant $C$ depends on $C_1,C_2,C_3,C_4, p$ and $\|u\|_{\rm BMO}$.
}
\end{lemma}
\noindent{\it Proof}~~
For any ball $B\subset\R^n$, we set $B_1=\{x\in B|\, |u(x)-u_B|<\frac{1}{2}u_B\}$ and $B_2=B\backslash B_1$. By the Nirenberg inequality \cite{Grafakos} we have
$|B_2|/|B|\leq e^{-Cu_B/\|u\|_{\rm BMO}}$. Clearly, $\rho(u)\geq C\rho(u_B)$ on $B_1$. Therefore,
\begin{align*}
\frac{1}{|B|}\int_{B}\rho(u)\d x
&\leq \frac{C}{|B|}\int_{B}(1+|u-u_B|^p)\d x + \frac{C}{|B|}\int_{B}|u_B|^p\d x\\
&\leq C + C|u_B|^p\leq C\rho(u_B) ,\\[8pt]
\frac{1}{|B|}\int_{B}\frac{1}{\rho(u)}\d x
&\leq \frac{1}{|B|}\int_{B_1}\frac{1}{\rho(u)}\d x + \frac{C|B_2|}{|B|}\\
&\leq \frac{C}{\rho(u_B)}+e^{-Cu_B/\|u\|_{\rm BMO}}
\leq \frac{C}{\rho(u_B)} .
\end{align*}
The last two inequalities imply that $\rho(u)$ is a $A_2$ weight.
\qed

The following lemma concerns maximal regularity of parabolic
equations, which is an application of the maximal regularity of \cite{Wood} and  \cite{JK} (with the perturbation method for the
treatment of operators with merely continuous coefficients).
\begin{lemma}\label{maxregpara}
{\it Let $u$ be the solution of the
parabolic problem {\rm (\ref{parabEq0})} in $\R^n$ $(n=2,3)$
with the Dirichlet boundary/initial conditions $u\equiv g\equiv f_0\equiv 0$, and assume that the coefficient matrix $A$ is continuous. Then we have
$$
\|u\|_{L^p(I;W^{1,q}(\Omega))}\leq
C_{p,q}\|\vec{f}\|_{L^p(I;L^q(\Omega))}
$$
for some $q>n$ and any $1<p<\infty$. The constant $C_{p,q}$ depends only on $p,q$, $K$, the domain $\Omega$ and the modulo of continuity of $A$.
}
\end{lemma}

The analogus result for elliptic equations is given below, which can be proved by applying the $W^{1,q}$ estimate of \cite{JK} with a perturbation argument.
\begin{lemma}\label{maxregellp}
{\it
Let $A_{ij}$, $i,j=1,\cdots,n$, be continuous functions defined on
$\Omega$, satisfying
\begin{align*}
K^{-1}|\xi|^2\leq \sum_{i,j=1}^n A_{ij}(x)\xi_i\xi_j \leq K|\xi|^2
,\quad\mbox{for all}~\xi\in\R^n,~~\mbox{a.e. $x\in\R^{n}$ $(n=2,3)$},
\end{align*}
where $K$ is a
positive constant.
Let $u$ be the solution of the
elliptic equation
$$
-\nabla\cdot(A\nabla
u)=\nabla\cdot\vec{f}\quad\mbox{in}~~\Omega,
$$
with the Dirichlet boundary/initial conditions $u=0$
on $\partial\Omega$. Then we have
$$
\|u\|_{W^{1,q}(\Omega)}\leq
C_{q}\|\vec{f}\|_{L^q(\Omega)}
$$
for some $q>n$. The constant $C_{q}$ depends only on $q$, $\Lambda$, the domain $\Omega$ and the modulo of continuity of $A$.
}
\end{lemma}

The following lemma is concerned with H\"{o}lder estimates for inhomogeneous parabolic equations \cite{AS}, which is also a consequence of Theorem \ref{Mainresult}.
\begin{lemma}\label{Mainresult4}
{\it The solution of {\rm(\ref{parabEq0})} with $u_0\equiv g\equiv 0$
satisfies that
\begin{align*}
&\|u\|_{C^{\alpha,\alpha/2}(\overline{\Omega}_T)}
\leq C(
\|f_0\|_{L^{p}((0,T);L^{q}(\Omega))}+\|{\vec
f}\|_{L^{2p}((0,T);L^{2q}(\Omega))}),
\end{align*}
for some $0<\alpha<1$, provided $1\leq p,q\leq\infty$ and $\,
2/p+n/q<2$. }
\end{lemma}

The following lemma concerns an estimate of $\nabla u$ in the Morrey space for the parabolic equation (\ref{parabEq0}), which was proved in \cite{Yin} for $u_0\equiv g\equiv f_0\equiv 0$.
\begin{lemma}\label{Mainrasdh8}
{\it The solution of {\rm(\ref{parabEq0})} with $f_0\equiv 0$
satisfies that
\begin{align*}
&\|\nabla u\|_{L^{2,n/(n+2)}_{\rm para}(\Omega_T)} \leq C(\|{\vec
f}\|_{L^{2,n/(n+2)}_{\rm para}(\Omega_T)} +\|\nabla
g\|_{L^{2,n/(n+2)}_{\rm
para}(\Omega_T)}+\|\partial_tg\|_{L^{2,n/(n+2)}_{\rm
para}(\Omega_T)}+\|u_0\|_{L^\infty(\Omega)}).
\end{align*}
}
\end{lemma}

\subsection{Construction of approximating solutions}

For the non-degenerate problem, the existence of a $C^\alpha$ solution was proved by Yuan and Lin \cite{YL,Yuan2}. Based on their result, for any given $\varepsilon>0$, there exists a weak solution $(u^\varepsilon,\phi^\varepsilon)$ such that $\phi^\varepsilon\in L^\infty((0,T);H^1(\Omega))$
and $u^\varepsilon\in C^{\alpha,\alpha/2}(\overline\Omega_T)\cap L^2((0,T);H^1(\Omega))$, to the following equations
\begin{align}\label{sdbgkfo99033}
\left\{
\begin{array}{ll}
\displaystyle
\frac{\partial u^\varepsilon}{\partial t}-\nabla\cdot(\kappa(u^\varepsilon)
\nabla
u^\varepsilon)=\nabla\cdot[
(\sigma(u^\varepsilon)+\varepsilon)\phi^\varepsilon\nabla\phi^\varepsilon]
&\mbox{in}~~\Omega,
\\[5pt]
u^\varepsilon=g  & \mbox{on}~~\partial\Omega ,\\
u^\varepsilon(x,0)=u_0(x)  & \mbox{for}~~x\in \Omega ,
\end{array}\right.
\end{align}
\begin{align}\label{sfdhyy733}
\left\{
\begin{array}{ll}
\displaystyle -\nabla\cdot\big(
(\sigma(u^\varepsilon)+\varepsilon)
\nabla\phi^\varepsilon\big)=0 &\mbox{in}~~\Omega,
\\[3pt]
\phi^\varepsilon=h & \mbox{on}~~\partial\Omega .
\end{array}\right.
\end{align}
We also note that, by the maximum principle, the solution $u^\varepsilon$ of  (\ref{sdbgkfo99033}) satisfies that
\begin{align}
u^\varepsilon\geq c:=
\min\big(\min_{x\in\Omega}u_0(x),\min_{x\in\partial\Omega}g(x)\big)>0
,
\end{align}
and the solution $\phi$ of (\ref{sfdhyy733}) satisfies that
\begin{align}\label{sdfjkljoo}
\|\phi^\varepsilon\|_{L^\infty(\Omega_T)}\leq \|h\|_{L^\infty(\Gamma_T)} .
\end{align}
By the hypotheses (H1)-(H2), we have
\begin{align}
\kappa_0\leq\kappa(u^\varepsilon)\leq \kappa_1,\qquad \varepsilon\leq\sigma(u^\varepsilon)+\varepsilon\leq 2\sigma_0:=\sup_{s\geq c}\,\sigma(s),
\end{align}
for some positive constants $\kappa_0,\kappa_1$ and $\sigma_0$, where we choose $\varepsilon<\sigma_0$.

\begin{proposition}\label{fdsjklhiy87979}
{\it
The solution $(u^\varepsilon,\phi^\varepsilon)$ of {\rm (\ref{sdbgkfo99033})-(\ref{sfdhyy733})} satisfies that
\begin{align*}
&\|u^\varepsilon\|_{C^{\alpha,\alpha/2}(\overline\Omega_T)}
+\|u^\varepsilon\|_{L^p((0,T);W^{1,q}(\Omega))}
+\|\partial_tu^\varepsilon\|_{L^p((0,T);W^{-1,q}(\Omega))}
+\|\phi^\varepsilon\|_{L^\infty((0,T);W^{1,q}(\Omega))}\leq C ,
\end{align*}
and
$$
\|\phi^\varepsilon\|_{C^{\alpha,\alpha/2}(\overline B_R\times[0,T])}\leq C_{{\rm dist}(\overline B_R,\partial\Omega)}
$$
for any closed ball $\overline B_R\subset\Omega$, where the constants $C$ and $C_{{\rm dist}(\overline B_R,\partial\Omega)} $ are independent of $\varepsilon$.
}
\end{proposition}

\noindent{\it Proof}~~~
First, we show that $\sigma(u^\varepsilon)+\varepsilon$ is a $A_2$ weight, uniformly with respect to time and $\varepsilon$.

Let $x_0\in\Omega$, $t_0>0$ and let
$R_0=\frac{1}{2}\min(\sqrt{t_0},{\rm dist}(x_0,\partial\Omega))$.
For any ball $B_R$ of radius $R$ centered at $x_0$, we let $\zeta$
be a smooth function defined on $\R^n$ which satisfies $0\leq
\zeta\leq 1$, $\zeta=1$ in $B_R$ and $\zeta=0$ outside $B_{2R}$. For
any interval $I_R=(t_0-R^2,t_0]$, we let $\chi$ be a smooth function
defined on $\R$ which satisfies $0\leq \chi\leq 1$, $\chi=1$ on
$I_R$ and $\chi=0$ on
$(-\infty,t_0-4R^2]$.
Let $Q_R=B_R\times I_R$ so that $(u^\varepsilon,\phi^\varepsilon)$ is a
solution of (\ref{sdbgkfo99033})-(\ref{sfdhyy733}) in $Q_{2R_0}$.
Multiplying (\ref{sfdhyy733}) by $\varphi=\phi^\varepsilon\zeta^2$, we obtain
\begin{align*}
\int_{B_R}(\sigma(u^\varepsilon)+\varepsilon)|\nabla\phi^\varepsilon|^2\d
x\leq\int_{B_{2R}}
(\sigma(u^\varepsilon)+\varepsilon)
|\phi^\varepsilon|^2|\nabla\zeta|^2\d x\leq
C\|\phi^\varepsilon\|_{L^\infty(\Omega)}^2R^{n-2} .
\end{align*}
Integrating the above inequality with respect to time and using (\ref{sdfjkljoo}), we get
\begin{align}\label{jkdfabe}
\iint_{Q_R}(\sigma(u^\varepsilon)+\varepsilon)|\nabla\phi^\varepsilon|^2\d x\d t\leq
C\|h\|_{L^\infty(\Gamma_T)}^2R^n.
\end{align}
Similarly, for $x_0\in\partial\Omega$, $t_0>0$, $R<\frac{1}{2}\min(\sqrt{t_0},R_\Omega)$, $B_R:=B_R(x_0)\cap\Omega$ and $Q_R=Q_R(x_0,t_0)\cap\Omega_T$, we also have (\ref{jkdfabe}).
From the last inequality we see that
\begin{align}\label{jkdfabe1}
\big\|\sqrt{\sigma(u^\varepsilon)+\varepsilon}
\nabla\phi^\varepsilon\big \|_{L^{2,n/(n+2)}_{\rm
para}(\Omega_T)}\leq C .
\end{align}

By Theorem \ref{Mainresult}, the solution of (\ref{sdbgkfo99033}) satisfies that
\begin{align}
\|u^\varepsilon\|_{L^\infty((0,T);\overline{{\rm BMO}})}\leq
C\big\|\sqrt{\sigma(u^\varepsilon)+\varepsilon}
\nabla\phi^\varepsilon\big \|_{L^{2,n/(n+2)}_{\rm
para}(\Omega_T)}+C\|u_0\|_{L^\infty(\Omega)}+C\|g\|_{L^\infty(\Gamma_T)}\leq
C .
\end{align}
Applying Lemma \ref{Mainrasdh8} to the equation (\ref{sdbgkfo99033}) and using (\ref{jkdfabe1}), we derive that
\begin{align}\label{sfdjkl79009}
&\|\nabla u^\varepsilon\|_{L^{2,n/(n+2)}_{\rm para}(\Omega_T)} \leq
C\big\|\sqrt{\sigma(u^\varepsilon)+\varepsilon}
\nabla\phi^\varepsilon\big \|_{L^{2,n/(n+2)}_{\rm
para}(\Omega_T)}+C\leq C .
\end{align}

We extend the function $u^\varepsilon$ defined on $\Omega$ to $\R^n$ by setting
$u^\varepsilon(x)=c$ for $x\in\R^n\backslash\Omega$ so that
$$\|u^\varepsilon\|_{L^\infty((0,T); {\rm BMO}(\R^n))}\leq C.$$

Since (\ref{dfsjljiowuo}) holds, from Lemma \ref{BMOandrecipr} we see that $\rho(u^\varepsilon)$ (and also $\sigma(u^\varepsilon)=1/\rho(u^\varepsilon)$) is a $A_2$ weight uniform with respect to time and $\varepsilon$.
It follows that, for any ball $B\subset \R^n$,
\begin{align*}
&\biggl(\frac{1}{|B|}\int_B(\sigma(u^\varepsilon)+\varepsilon) \d x\biggl)\biggl(\frac{1}{|B|}\int_B\frac{1}{\sigma(u^\varepsilon)+\varepsilon}\d x\biggl)  \\\
&=\biggl(\frac{1}{|B|}\int_B\sigma(u^\varepsilon) \d x\biggl)\biggl(\frac{1}{|B|}\int_B\frac{1}{\sigma(u^\varepsilon)+\varepsilon}\d x\biggl)+\frac{1}{|B|}\int_B\frac{\varepsilon}{\sigma(u^\varepsilon)+\varepsilon}\d x\\
&\leq \biggl(\frac{1}{|B|}\int_B\sigma(u^\varepsilon) \d x\biggl)\biggl(\frac{1}{|B|}\int_B\frac{1}{\sigma(u^\varepsilon)}\d x\biggl)+1\\
&\leq C ,
\end{align*}
which says that $\sigma(u^\varepsilon)+\varepsilon$ is also a $A_2$ weight, uniform with respect to time and $\varepsilon$.

Secondly, we estimate the H\"{o}lder norms of $\phi^\varepsilon$ and $u^\varepsilon$, respectively. In fact,
from \cite{FKS} we know that any solution of
the elliptic equation (\ref{sfdhyy733}) with the $A_2$ coefficient
$\sigma(u^\varepsilon)+\varepsilon$ satisfies the H\"{o}lder estimates:
\begin{align}\label{holderest}
\|\phi^\varepsilon(\cdot,t)\|_{C^\alpha(\overline\Omega)}\leq
C\|h^\varepsilon(\cdot,t)\|_{C^\alpha(\partial\Omega)}
\leq C,
\quad\mbox{for}~~ t\in(0,T),\quad\forall~\alpha\in(0,\alpha_0),
\end{align}
for some fixed constant $\alpha_0\in (0,1)$.

We proceed to the H\"{o}lder estimate of $u^\varepsilon$. For any fixed $x_0\in\Omega$,
we decompose the function $u^\varepsilon$ as
$u^\varepsilon=u_1^\varepsilon+u_2^\varepsilon$, where $u_1^\varepsilon$ and $u_2^\varepsilon$ are weak solutions of the equations
\begin{align*}
\left\{
\begin{array}{ll}
\displaystyle
\frac{\partial u_1^\varepsilon}{\partial t}-\nabla\cdot(\kappa(u^\varepsilon)
\nabla
u_1^\varepsilon)=0
&\mbox{in}~~\Omega,
\\[5pt]
u_1^\varepsilon=g & \mbox{on}~~\partial\Omega ,\\
u_1^\varepsilon(x,0)=u_0(x) & \mbox{for}~~x\in \Omega ,
\end{array}\right.
\end{align*}
and
\begin{align*}
\left\{
\begin{array}{ll}
\displaystyle
\frac{\partial u_2^\varepsilon}{\partial t}-\nabla\cdot(\kappa(u^\varepsilon)
\nabla
u_2^\varepsilon)=\nabla\cdot\big[
(\phi^\varepsilon-\phi^\varepsilon(x_0,t))
(\sigma(u^\varepsilon)+\varepsilon)
\nabla\phi^\varepsilon\big]
&\mbox{in}~~\Omega,
\\[5pt]
u_2^\varepsilon=0 & \mbox{on}~~\partial\Omega ,\\
u_2^\varepsilon(x,0)=0 & \mbox{for}~~x\in \Omega ,
\end{array}\right.
\end{align*}
respectively.
By the De Giorgi--Nash--Moser estimates, we have
$$\|u_1^\varepsilon\|_{C^{\alpha,\alpha/2}(\overline\Omega_T)}\leq C(\|g\|_{C^{\alpha,\alpha/2}(\overline\Gamma_T)}+\|u_0\|_{C^\alpha(\overline\Omega)}).$$
and in order to estimate $\|u_2^\varepsilon\|_{C^{\alpha,\alpha/2}(\overline\Omega_T)}$, we set $\vec{f}=(\phi^\varepsilon-\phi^\varepsilon(x_0,t))
(\sigma(u^\varepsilon)+\varepsilon)
\nabla\phi^\varepsilon$ and apply Proposition \ref{dnflkqaghu}--Proposition \ref{dnflkqsady63}. We see that
for $x_0\in\Omega$, $t_0>0$, $0<2\rho\leq R\leq\min\big({\rm
dist}(x_0,\partial\Omega),\sqrt{t_0}\big)$,
we have
\begin{align*}
&\frac{1}{\rho^{n+2+2\alpha}}
\|u_2^\varepsilon-(u_2^\varepsilon)_{Q_\rho}\|_{L^2(Q_\rho)}^2\\
&\leq
C\biggl(\frac{1}{R^{n+2+2\alpha}}
\|u_2^\varepsilon-\theta\|_{L^2(Q_{R})}^2
+\frac{1}{R^{n+2\alpha}}\|\vec{f}\|_{L^2(Q_R)}^2
\biggl), \\
&\leq
C\biggl(\frac{1}{R^{n+2+2\alpha}}
\|u_2^\varepsilon-\theta\|_{L^2(Q_{R})}^2
+\frac{1}{R^{n}}\|\phi^\varepsilon
\|_{L^\infty(I;C^\alpha(\overline\Omega))}^2\big\|
\sqrt{\sigma(u^\varepsilon)+\varepsilon}
\nabla\phi^\varepsilon\big\|_{L^2(Q_R)}^2
\biggl) \\
&\leq
C\biggl(\frac{1}{R^{n+2+2\alpha}}
\|u_2^\varepsilon-\theta\|_{L^2(Q_{R})}^2
+\frac{1}{R^{n}}\big\|
\sqrt{\sigma(u^\varepsilon)+\varepsilon}
\nabla\phi^\varepsilon\big\|_{L^2(Q_R)}^2
\biggl) .
\end{align*}
Similarly, for $x_0\in\Omega$, $t_0=0$, $Q_R=B_R(x_0)\times[0,R^2]$ and $0<\rho<R\leq\min\big({\rm
dist}(x_0,\partial\Omega),\sqrt{T}\big)$, we
have
\begin{align*}
&\frac{1}{\rho^{n+2+2\alpha}}
\|u_2^\varepsilon\|_{L^2(Q_\rho)}^2\leq C\biggl(\frac{1}{R^{n+2+2\alpha}}
\|u_2^\varepsilon\|_{L^2(Q_{R})}^2
+\frac{1}{R^{n}}\big\|
\sqrt{\sigma(u^\varepsilon)+\varepsilon}
\nabla\phi^\varepsilon\big\|_{L^2(Q_R)}^2
\biggl)
. \nn
\end{align*}
For $x_0\in\partial \Omega$, $t_0>0$, $Q_R=B_R(x_0)\cap\Omega\times(t_0-R^2,t_0]$ and $0<\rho<R\leq{\rm
dist}(R_\Omega,\sqrt{t_0})$, we have
\begin{align*}
&\frac{1}{\rho^{n+2+2\alpha}}
\|u_2^\varepsilon\|_{L^2(Q_\rho)}^2\leq
C\biggl(\frac{1}{R^{n+2+2\alpha}}\|u_2^\varepsilon\|_{L^2(Q_{R})}^2
+\frac{1}{R^{n}}\big\|
\sqrt{\sigma(u^\varepsilon)+\varepsilon}
\nabla\phi^\varepsilon\big\|_{L^2(Q_R)}^2
\biggl) . \nn
\end{align*}
For $x_0\in\partial \Omega$, $t_0=0$, $Q_R=B_R(x_0)\cap\Omega\times[0,R^2]$ and $0<\rho<R\leq
\min(R_\Omega,\sqrt{T})$, we have
\begin{align*}
&\frac{1}{\rho^{n+2+2\alpha}}
\|u_2^\varepsilon\|_{L^2(Q_\rho)}^2\leq
C\biggl(\frac{1}{R^{n+2+2\alpha}}
\|u_2^\varepsilon\|_{L^2(Q_{R})}^2
+\frac{1}{R^{n}}
\big\|
\sqrt{\sigma(u^\varepsilon)+\varepsilon}
\nabla\phi^\varepsilon\big\|_{L^2(Q_R)}^2
\biggl)
. \nn
\end{align*}

Combining the last four inequalities and following the outline of Section \ref{glbest1}, we can derive that
\begin{align*}
\|u_2^\varepsilon\|_{{\cal L}^{2,1+2\alpha/(n+2)}_{\rm
para}(\Omega_T)} \leq C\big\|
\sqrt{\sigma(u^\varepsilon)+\varepsilon} \nabla\phi^\varepsilon\big
\|_{L^{2,n/(n+2)}_{\rm para}(Q_R)}.
\end{align*}
With (\ref{jkdfabe1}) and the equivalence relation ${\cal L}^{2,1+2\alpha/(n+2)}_{\rm para}(\Omega_T)\cong C^{\alpha,\alpha/2}(\overline\Omega_T)$, we see that
\begin{align*}
\|u_2^\varepsilon\|_{C^{\alpha,\alpha/2}
(\overline\Omega_T)} \leq
C .
\end{align*}
Therefore,
\begin{align}
\|u^\varepsilon\|_{C^{\alpha,\alpha/2}
(\overline\Omega_T)}\leq \|u_1^\varepsilon\|_{C^{\alpha,\alpha/2}
(\overline\Omega_T)}+\|u_2^\varepsilon\|_{C^{\alpha,\alpha/2}
(\overline\Omega_T)} \leq
C .
\end{align}

Thirdly, we present $W^{1,q}$ estimates of $\phi^\varepsilon$ and $u^\varepsilon$.
Note that the last inequality implies that
\begin{align}
C^{-1}\leq \sigma(u^\varepsilon)+\varepsilon\leq C,\quad \|\sigma(u^\varepsilon)+\varepsilon
\|_{C^{\alpha,\alpha/2}(\overline\Omega_T)}
\leq C,\quad  \|\kappa(u^\varepsilon)
\|_{C^{\alpha,\alpha/2}(\overline\Omega_T)}
\leq C.
\end{align}
With the H\"{o}lder estimates of  $\sigma(u^\varepsilon)+\varepsilon$ and $\kappa(u^\varepsilon)$, we
apply Lemma \ref{maxregpara} -- \ref{maxregellp} and derive that
\begin{align}
&\|\phi^\varepsilon\|_{L^\infty((0,T);W^{1,q}(\Omega))} \leq
C\|h\|_{L^\infty((0,T);W^{1,q}
(\Omega))}\leq C
,\\
&\|u^\varepsilon\|_{L^p((0,T);W^{1,q}(\Omega))} \leq
C_p\|\phi^\varepsilon\|_{L^p((0,T);W^{1,q}(\Omega))}+C_p
\leq C_p ,
\end{align}
for some $q>n$ and any $1< p<\infty$. From the equation (\ref{sdbgkfo99033}) we also see that
\begin{align}
\|\partial_tu^\varepsilon\|_{L^p(I;W^{-1,q}(\Omega))} \leq
C(\|u^\varepsilon\|_{L^p((0,T);W^{1,q}(\Omega))}+\|\nabla\phi^\varepsilon\|_{L^p((0,T);W^{1,q}(\Omega))})\leq C.
\end{align}

Finally, we estimate the interior space-time H\"{o}lder norm of $\phi^\varepsilon$, which is used to obtain pointwise convergence of the approximating solutions in the next subsection.
For the simplicity of notations, we set $A^\varepsilon=\sigma(u^\varepsilon)+\varepsilon$. From (\ref{sfdhyy733}) we see that
\begin{align*}
\displaystyle -\nabla\cdot\biggl(
A^\varepsilon(x,t_1)
\nabla[\phi^\varepsilon(x,t_1)
-\phi^\varepsilon(x,t_2)]\biggl)
=\nabla\cdot\biggl(
(A^\varepsilon(x,t_1)-A^\varepsilon(x,t_2))
\nabla\phi^\varepsilon(x,t_2)\biggl) .
\end{align*}
By applying the interior $W^{1,q}$ estimate to the above equation, we find that for any closed ball $\overline B_R$ contained in $\Omega$ there holds
\begin{align*}
\|\phi^\varepsilon(x,t_1)
-\phi^\varepsilon(x,t_2)\|_{
L^\infty((0,T);W^{1,q}(B_R))} &\leq
C_{{\rm dist}(\overline B_R,\partial\Omega)}
\|A^\varepsilon(x,t_1)-A^\varepsilon(x,t_2)
\|_{L^\infty(\Omega_T)}\\
&\leq
C_{{\rm dist}(\overline B_R,\partial\Omega)}
\|A^\varepsilon\|_{C^{\alpha,\alpha/2}
(\overline\Omega_T)}|t_1-t_2|^{\alpha/2} ,
\end{align*}
which reduces to
\begin{align*}
\|\phi^\varepsilon\|_{
C^{\alpha/2}([0,T]
;W^{1,q}(B_R))} \leq
C_{{\rm dist}(\overline B_R,\partial\Omega)} .
\end{align*}
Since $W^{1,q}(B_R)\hookrightarrow C^\alpha(\overline\Omega)$, the last inequality implie that
\begin{align}
\|\phi^\varepsilon\|_{
C^{\alpha,\alpha/2}(\overline B_R\times[0,T])} \leq
C_{{\rm dist}(\overline B_R,\partial\Omega)} .
\end{align}

The proof of Proposition \ref{fdsjklhiy87979} is complete.\qed

\subsection{Existence of solution}
Since $C^{\alpha,\alpha/2}
(\overline\Omega_T)$ is compactly embedded into $C(\overline\Omega_T)$ and $C^{\alpha,\alpha/2}(\overline B_R\times[0,T])$ is compactly embedded into $C(\overline B_R\times[0,T])$
there exist functions $u\in C^{\alpha,\alpha/2}
(\overline\Omega_T)$, $\phi\in L^\infty(I;W^{1,q}(\Omega))$ with $\phi\in C^{\alpha,\alpha/2}(\overline B_R\times[0,T])$ for any closed ball $\overline B_R$ contained in $\Omega$, and a sequence $\varepsilon_k\rightarrow 0$, such that
$u^{\varepsilon_k}$ converges to $u$ in the norm of $C(\overline\Omega_T)$, $u^{\varepsilon_k}$ converges weakly to $u$ in $L^p(I;W^{1,q}(\Omega))$, $\partial_tu^{\varepsilon_k}$ converges weakly to $\partial_tu$ in $L^p(I;W^{-1,q}(\Omega))$, $\phi^{\varepsilon_k}$ converges weakly$^*$ to $\phi$ in $L^\infty(I;W^{1,q}(\Omega))$, and $\phi^{\varepsilon_k}$ converges to $\phi$ pointwise uniformly in each compact subset of $\Omega\times[0,T]$.

From (\ref{sfdhyy733}) we see that
\begin{align*}
\int_\Omega(\sigma(u^{\varepsilon_k})+\varepsilon_k)
\nabla\phi^{\varepsilon_k}\cdot\nabla\varphi \,\d x=0
\quad\mbox{for any $\varphi\in H^1_0(\Omega)$} .
\end{align*}
By taking the limit $k\rightarrow\infty$, we obtain
\begin{align}\label{sdfhjkyo79}
\int_\Omega\sigma(u)\nabla\phi\cdot\nabla\varphi \,\d x=0, \quad\mbox{for any $\varphi\in H^1_0(\Omega)$ and a.e. $t\in (0,T)$.}
\end{align}
Therefore, for any function $v\in C^\infty_0(\Omega)$,
\begin{align*}
\lim_{k\rightarrow\infty}
\int_\Omega\nabla\cdot\biggl(
\phi^{\varepsilon_k}(\sigma(u^{\varepsilon_k})+\varepsilon_k)
\nabla\phi^{\varepsilon_k}\biggl)v \d x
&=
-\lim_{k\rightarrow\infty}\int_\Omega
\phi^{\varepsilon_k}(\sigma(u^{\varepsilon_k})+\varepsilon_k)
\nabla\phi^{\varepsilon_k}\cdot\nabla v \d x\\
&=
-\int_\Omega
\phi\sigma(u)
\nabla\phi\cdot\nabla v\d x\\
&=
-\int_\Omega
\sigma(u)
\nabla\phi\cdot[\nabla(\phi v)-v\nabla\phi] \d x\\
&=\int_\Omega
\sigma(u)
|\nabla\phi|^2v\d x .
\end{align*}

From (\ref{sdbgkfo99033}) we know that for any $v\in L^\infty((0,T);C^\infty_0(\Omega))$,
\begin{align*}
\int_0^T\int_\Omega\frac{\partial u^{\varepsilon_k}}{\partial t}v \d x\d t+\int_0^T\int_\Omega
\kappa(u^{\varepsilon_k})_{\varepsilon_k}
\nabla
u^\varepsilon\cdot\nabla v\d x\d t\\
=\int_0^T\int_\Omega\nabla\cdot\biggl(
\phi^{\varepsilon_k}\frac{1}{
\rho(u^{\varepsilon_k})_{\varepsilon_k}}
\nabla\phi^{\varepsilon_k}\biggl)v\,\d x\d t .
\end{align*}
By taking the limit $k\rightarrow\infty$, we get
\begin{align}\label{sdfhjkyoasdf}
\int_0^T\int_\Omega\frac{\partial u}{\partial t}v\, \d x\d t+\int_0^T\int_\Omega
\kappa(u)
\nabla
u\cdot\nabla v\, \d x\d t=\int_0^T\int_\Omega\sigma(u)
|\nabla\phi|^2v\,\d x\d t .
\end{align}

From the regularity of $u$ and $\phi$, we know that the equations (\ref{sdfhjkyo79})-(\ref{sdfhjkyoasdf}) actually hold for any $\varphi\in H^1_0(\Omega)$ and $v\in L^2((0,T);H^1_0(\Omega))$.

To conclude, we have proved the existence of a weak solution $(u,\phi)$ to the equations (\ref{e-heat-1})-(\ref{BC}) with the regularity (\ref{sdfjhkl79031}).

\subsection{Uniqueness of solution}

Suppose that $(u_1,\phi_1)$ and $(u_2,\phi_2)$ are two pairs of
solutions to the initial-boundary value problem
(\ref{e-heat-1})-(\ref{BC}), both satisfying  (\ref{sdfjhkl79031}). Let $\bar u=u_1-u_2$ and $\bar
\phi=\phi_1-\phi_2$. Then $\bar u$ and $\bar \phi$ are weak solutions to the equations
\begin{align}
&\frac{\partial \bar u}{\partial
t}-\nabla\cdot(\kappa(u_1)\nabla\bar
u)=\nabla\cdot((\kappa(u_1)-\kappa(u_2))\nabla
u_2) \nonumber \\
&~~~+(\sigma(u_1)-\sigma(u_2))|\nabla\phi_1|^2+\sigma(u_2)\nabla(\phi_1+\phi_2)\cdot\nabla\bar\phi
\label{diffEqu}\\[5pt]
&-\nabla\cdot(\sigma(u_1)\nabla\bar\phi)
=\nabla\cdot\Big((\sigma(u_1)-\sigma(u_2))\nabla\phi_2\Big),
\label{diffEqphi}
\end{align}
with the following boundary and initial conditions:
\begin{align}
\begin{array}{ll}
\bar u(x,t)=0,\quad \bar\phi(x,t)=0~~
&\mbox{for}~~x\in\partial\Omega,~~t\in[0,T],\\[3pt]
\bar u(x,0)=0~~ &\mbox{for}~~x\in\Omega .
\end{array}
\end{align}

For any $\tau\in(0,T)$, we denote $I_\tau=(0,\tau)$ and
$\Omega_\tau=\Omega\times I_\tau$. By applying Lemma \ref{Mainresult4}
to the parabolic equation (\ref{diffEqu}), we see that for $q>n$ there
exists $1<p<\infty$ such that
\begin{align*}
\|\bar u\|_{L^\infty(\Omega_\tau)}&\leq
C\|(\kappa(u_1)-\kappa(u_2))\nabla u_2\|_{L^p(I_\tau;L^q(\Omega))} \nn\\
&~~~+C\|(\sigma(u_1)-\sigma(u_2))|\nabla\phi_1|^2\|_{L^p(I_\tau;L^{q/2}(\Omega))} \nn\\
&~~~+C\|\sigma(u_2)\nabla(\phi_1+\phi_2)\cdot\nabla\bar\phi\|_{L^p(I_\tau;L^{q/2}(\Omega))} \nn\\
&\leq C\|\bar u\|_{L^\infty(\Omega_\tau)}\big(\tau^{1/2p}\|\nabla
u_2\|_{L^{2p}(I_\tau;L^q(\Omega))}+\|\nabla\phi_1\|_{L^{2p}(I_\tau;L^q(\Omega))}^2\big) \nn\\
&~~~+C\tau^{1/2p}\|\nabla(\phi_1+\phi_2)\|_{L^\infty(I_\tau;L^q(\Omega))}
\|\nabla\bar\phi\|_{L^{2p}(I_\tau;L^q(\Omega))} \nn\\
&\leq C\tau^{1/2p}\|\bar
u\|_{L^\infty(\Omega_\tau)}+C\tau^{1/2p}\|\nabla\bar\phi\|_{L^\infty(I_\tau;L^q(\Omega))}
,
\end{align*}
where the constant $C$ is independent of $\tau$.
With the H\"{o}lder regularity of $u_1$, by applying the $W^{1,q}$
estimates to (\ref{diffEqphi}), we obtain
\begin{align*}
\|\nabla\bar\phi\|_{L^\infty(I_\tau;L^q(\Omega))}&\leq
C\|(\sigma(u_1)-\sigma(u_2))\nabla\phi_2\|_{L^\infty(I_\tau;L^q(\Omega))}
\leq C\|\bar u\|_{L^\infty(\Omega_\tau)} .
\end{align*}
There exists $T_0$ such that for $\tau<T_0$, the last two
inequalities imply that
\begin{align*}
\|\bar
u\|_{L^\infty(\Omega_\tau)}+\|\nabla\bar\phi\|_{L^\infty(I_\tau;L^q(\Omega))}
=0 .
\end{align*}
By dividing the interval $(0,T)$ into small parts $(T_k,T_{k+1}]$,
$k=0,1,\cdots$, each part satisfying $T_{k+1}-T_k<T_0$, we find that
$\bar u(\cdot,T_k)\equiv\bar \phi(\cdot,T_k)\equiv0$ implies that
$\bar u(\cdot,t)\equiv\bar \phi(\cdot,t)\equiv0$ for
$t\in[T_k,T_{k+1}]$. This proves the uniqueness of solution.

\medskip

\section{Conclusions}
\setcounter{equation}{0}

In this paper, we proved global existence and uniqueness of a weak
solution to the degenerate thermistor problem by establishing a
uniform-in-time BMO estimate for parabolic equations with possibly
discontinuous coefficients. The physical hypothesis (H1)-(H2) are
satisfied by metals and some semiconductors. The BMO estimate of
parabolic equations established in this paper may be applied to many
other equations of mathematical physics.

%

\end{document}